\documentclass[12pt]{article}
\usepackage{mathrsfs}
\usepackage{amssymb}
\usepackage{amsfonts}

\usepackage{amsmath}
\usepackage{latexsym,amsfonts,amssymb}

\newcommand{\meas}{{\rm meas}}
\newcommand{\epsln}{\varepsilon}
\newcommand{\sgn}{{\rm \,sgn\,}}
\newcommand{\supp}{{\rm supp}}

\renewcommand{\theequation}{\thesection.\arabic{equation}}

\title{Weak Continuity of the Flow Map for the Benjamin-Ono Equation on the Line}
\author{Shangbin Cui$^\dag$ \ \ and \ \ Carlos E. Kenig$^\ddag$\\ [0.1cm]
 {\small $^\dag$ Department of Mathematics, Sun Yat-Sen University,
  Guangzhou, Guangdong 510275}\\
 {\small People's Republic of China. E-mail: cuisb@yahoo.com.cn}\\ [0.1cm]
 {\small $^\ddag$ Department of Mathematics, University of Chicago,
  Chicago, IL 60637, USA.}\\
 {\small E-mail: cek@math.uchicago.edu}}
\date{}

\begin{document}
\maketitle

\begin{abstract}
  In this paper we show that the flow map of the Benjamin-Ono equation on the
  line is weakly continuous in $L^2(\mathbb{R})$, using ``local smoothing''
  estimates. $L^2(\mathbb{R})$ is believed to be a borderline space for the
  local well-posedness theory of this equation. In the periodic case, Molinet
  \cite{M1} has recently proved that the flow map of the Benjamin-Ono equation
  is not weakly continuous in $L^2(\mathbb{T})$. Our results are in line with
  previous work on the cubic nonlinear Schr\"{o}dinger equation, where
  Goubet and Molinet \cite{GoM} showed weak continuity in $L^2(\mathbb{R})$
  and Molinet \cite{M2} showed lack of weak continuity in $L^2(\mathbb{T})$.
\end{abstract}

\section{Introduction}
\setcounter{equation}{0}

  In this paper we study the weak continuity of the solution operator of the
  initial value problem for the Benjamin-Ono equation:
$$
\left\{
\begin{array}{l}
  \partial_t u+\mathcal{H}\partial_x^2 u+
  u\partial_x u=0, \quad
  x\in\mathbb{R},\;\;\; t\in\mathbb{R},\\
  u(x,0)=\phi(x), \quad x\in\mathbb{R},
\end{array}
\right.
\eqno{(1.1)}
$$
  where $\mathcal{H}$ represents the Hilbert transformation.

  The Benjamin-Ono equation (1.1) is a model for one-dimensional long waves in
  deep water (cf. \cite{B} and \cite{Ono}) and is completely integrable.
  Well-posedness of the problem (1.1) has been extensively studied by many
  authors, cf. \cite{BiL}, \cite{BP}, \cite{CKS}, \cite{GV1}--\cite{GV3},
  \cite{IK}, \cite{I}, \cite{KPV3}, \cite{KK}, \cite{KoTz}, \cite{M1}, \cite{M2},
  \cite{M4}, \cite{P}, \cite{T}, and the references therein. In particular, in
  \cite{IK}, it was proved that this problem is globally well-posed in
  $L^2(\mathbb{R})$. Thus, for any given $T>0$ there exists a mapping $S:
  L^2(\mathbb{R})\to C([-T,T],L^2(\mathbb{R}))$, which is Lipschitz continuous
  when restricted in any bounded sets in $L^2(\mathbb{R})$, such
  that for any $\phi\in L^2(\mathbb{R})$, the function $u(\cdot,t)=(S\phi)(t)=:
  S(t)\phi$ is a solution of the problem (1.1) in the time interval $[-T,T]$.
  In this paper we study the following problem: Is the operator $S(t):
  L^2(\mathbb{R})\to L^2(\mathbb{R})$ weakly continuous (for fixed $t$)?
  Note that since $S(t)$ is a nonlinear operator, we cannot give this question
  a positive answer by merely using the continuity of $S(t)$ in norm.

  Our motivation to study the above problem is inspired by the important
  series of works of Martel and Merle \cite{MM1}--\cite{MM4},
  which studied finite time blow-up and asymptotic stability and instability
  of solitary waves for the generalized Korteweg-de Vries equations, in
  critical and subcritical cases. One key step in their strategy in these
  works is a reduction to a nonlinear Liouville type theorem. Martel-Merle
  then reduce this nonlinear Liouville theorem to a corresponding linear
  one, involving the linearized operator around the solitary wave. It is in
  both these steps that the weak continuity of the flow map for generalized
  KdV in suitable Sobolev spaces plays a central role. Recently, by using
  a similar strategy, Kenig and Martel \cite{KM} established the asymptotic
  stability of solitons for the Benjamin-Ono equation (1.1) in the energy
  space $H^{1/2}(\mathbb{R})$. Thus, the weak continuity of the flow map in
  the energy space for the equation (1.1) is needed and it is established
  by these authors. The proof is very simple and reduces matters to the uniform
  continuity of the flow map for the Benjamin-Ono equation for data whose
  small frequencies coincide in a Sobolev space of strictly smaller index than
  $H^{1/2}(\mathbb{R})$, which depends on local well-posedness of the initial
  value problem (1.1) in $L^2(\mathbb{R})$ proved in the above mentioned work
  of Ionescu and Kenig \cite{IK}. Naturally, it would be desirable to prove the
  asymptotic stability of solitons for the Benjamin-Ono equation in
  $L^2(\mathbb{R})$. However, since no local well-posedness theory for this
  equation is available in Sobolev spaces of negative indices and it is strongly
  suspected that, in fact, uniform continuity of the flow map even restricted
  to data whose small frequencies coincide, must fail for Sobolev spaces of
  negative indices, the approach used in \cite{KM} does not work in the
  $L^2(\mathbb{R})$.

  Another interesting result which motivates this study is a recent work of
  Molinet \cite{M4}, in which the periodic initial-boundary value problem of
  the Benjamin-Ono equation was studied, and it was proved that the flow map
  of the periodic initial-boundary value problem of the Benjamin-Ono equation
  is not weakly continuous in $L^2(\mathbb{T})$, despite that such a problem
  is globally well-posed in $L^2(\mathbb{T})$, by another work of Molinet
  \cite{M2}.

  We would also like to mention a recent work of Goubet and Molinet
  \cite{GoM}, where a similar problem for the cubic nonlinear Schr\"{o}dinger
  equation on the line was studied. For this equation the global well-posedness in
  $L^2(\mathbb{R})$ was established in \cite{Ts}, while in \cite{KPV4}
  (focusing case) and \cite{CCT} (defocusing case) it was shown that the flow
  map is not uniformly continuous in any Sobolev space of negative index.
  Thus, the weak continuity in $L^2(\mathbb{R})$ of the flow map cannot be
  treated by the approach used in the works of Martel and Merle
  \cite{MM1}--\cite{MM4} and Kenig and Martel \cite{KM}. Goubet and Molinet
  \cite{GoM} affirmatively settled this problem by taking advantage of the
  ``local smoothing'' effect estimates together with a suitable uniqueness result.

  In this paper we establish the weak continuity in $L^2(\mathbb{R})$ of the
  flow map for the Benjamin-Ono equation (1.1). The main idea of the proof of this
  result is similar to that used in \cite{GoM}, i.e., we shall prove that the
  desired weak continuity is ensured by certain local compactness results
  coupled with suitable uniqueness. However, unlike the cubic nonlinear
  Schr\"{o}dinger case where local compactness is obtained from ``local
  smoothing'' effect estimates of the equation, in the present Benjamin-Ono
  case this will be derived from the properties of general functions in the
  space ${F}^\sigma$ in which local solutions of the problem (1.1) are
  constructed. Another interesting difference lies in the fact that, unlike
  the cubic nonlinear Schr\"{o}dinger case, the uniqueness for (1.1) is only
  established in \cite{IK} for limits of smooth solutions.

  To state our main result, we recall that $(p,q)$ is called an {\em admissible
  pair} for the operator $\partial_t+\mathcal{H}\partial_x^2$ if it satisfies
  the following conditions: $2\leq p\leq\infty$, $4\leq q\leq\infty$, and
  $2/q=1/2-1/p$. The main result of this paper reads as follows:
\medskip

  {\bf Theorem 1.1}\ \ {\em Assume that $\phi_n$ weakly converges to $\phi$ in
  $L^2(\mathbb{R})$. Let $u_n$ and $u$ be the solutions of the problem $(1.1)$
  with initial data $\phi_n$ and $\phi$, respectively, i.e., $u_n(\cdot,t)=S(t)
  \phi_n$ and $u(\cdot,t)=S(t)\phi$. Then given $T>0$, we have the following
  assertions:

  $(i)$ For any admissible pair $(q,p)$, $u_n$ weakly converges to $u$ in
  $L^q([-T,T],L^p(\mathbb{R}))$ $($in case either $q=\infty$ or $p=\infty$,
  weak convergence here refers to $*$-weak convergence$)$.

  $(ii)$ For any $|t|\leq T$, $u_n(t)$ weakly converges to $u(t)$ in
  $L^2(\mathbb{R})$. Moreover, this weak convergence is uniform for $|t|
  \leq T$ in the following sense: For any $\varphi\in L^2(\mathbb{R})$ we
  have
$$
  \lim_{n\to\infty}\sup_{|t|\leq T}|(u_n(t)-u(t),\varphi)|=0,
\eqno{(1.2)}
$$
  where $(\cdot,\;\cdot)$ denotes the inner product in $L^2(\mathbb{R})$.}
\medskip

  The arrangement of this paper is as follows: In Section 2 we give a
  review of the well-posedness result established in \cite{IK} and introduce
  the spaces used in the proof of this well-posedness result. In Section 3
  we derive some preliminary estimates. The proof of Theorem 1.1
  will be given in Section 4 after these preparations.

  Finally, we would like to give a remark on the modified Benjamin-Ono
  equation:
$$
  \partial_t u+\mathcal{H}\partial_x^2 u+
  u^2\partial_x u=0, \quad x\in\mathbb{R},\;\;\; t\in\mathbb{R}.
\eqno{(1.3)}
$$
  For this equation, it has been proved by Kenig and Takaoka in \cite{KT} that
  its initial value problem is globally well-posed in the Sobolev space
  $H^{1/2}(\mathbb{R})$, whereas the solution operator of a such problem is not
  uniformly continuous in any Sobolev space $H^{s}(\mathbb{R})$ of index
  $s<1/2$ (so that $H^{1/2}(\mathbb{R})$ is a borderline space for the local
  well-posedness theory of this equation). It is thus natural to ask if the
  flow map of this equation in $H^{1/2}(\mathbb{R})$ is weakly continuous. The
  answer to this question is affirmative and its proof is relatively easier,
  due to a priori regularities possessed by functions in the space $C([-T,T],
  H^{1/2}(\mathbb{R}))$. See the remark at the end of the paper.
\medskip

  {\bf Acknowledgement }\ \ This work on the part of the first author is
  financially supported by the National Natural Science Foundation under the
  grant number 10771223 as well as a fund from the Sun Yat-Sen University, and
  was performed when he was visiting the University of Chicago under financial
  support of China Scholarship Council. He would like to express his thanks to
  the Department of Mathematics of the University of Chicago for its hospitality
  during his visit. The second author is supported in part by NSF grant
  DMS-0456583.

  This second version of the manuscript (note the new title) deletes the section
  on the weak continuity of the flow map of the cubic nonlinear Schrodinger
  equation $in L^2(R)$ and gives a simpler proof of the weak continuity of the
  flow map of the modified Benjamin-Ono equation in $H^1/2(R)$. The authors
  are very indebted to an anonymous referee who pointed out to them that the work
  of Goubet and Molinet [11] already contained a proof for the cubic nonlinear
  Schrodinger and that a (simplified) version of the Goubet-Molinet proof could
  be used to give a very short proof of our mBO weak continuity result, which is
  presented in the last remarks of this paper.The authors would also like to thank
  Professor L.Molinet for helpful correspondence on these issues.
\medskip

{\bf Notations}:\ \
\begin{itemize}
\item For $1\leq p\leq\infty$, $\|\cdot\|_p$ denotes the norm in the
  Lebesque space $L^p(\mathbb{R})$.
\item For $1\leq p\leq\infty$, $1\leq q\leq\infty$ and a function $u=u(x,t)$
  ($x\in\mathbb{R}$, $t\in\mathbb{R}$), $\|u\|_{L_t^q L_x^p}$ and
  $\|u\|_{L_x^p L_t^q}$ denote norms of the mappings $t\to u(\cdot,t)$ and
  $x\to u(x,\cdot)$ in the spaces $L^q(\mathbb{R}_t,L^p(\mathbb{R}_x))$ and
  $L^p(\mathbb{R}_x,L^q(\mathbb{R}_t))$, respectively. In case that
  $\mathbb{R}_t$ is replaced by $[-T,T]$ or $\mathbb{R}_x$ by $[-R,R]$, the
  corresponding notation $L_t^q$ or $L_x^p$ in the norm notation will be
  replaced by $L_T^q$ or $L_R^p$, respectively, so that $\|u\|_{L_T^q L_R^p}$
  denotes the norm of the mapping $t\to u(\cdot,t)$ in the space
  $L^q([-T,T],L^p[-R,R])$, etc..
\item $\mathcal{F}$, $\mathcal{F}_1$ and $\mathcal{F}_2$ denote Fourier
  transformations in the varaibles $(x,t)$, $x$ and $t$, respectively; they
  will also be denoted as $\widetilde{\ \ }$, $\widehat{\ \ }^1$ and
  $\widehat{\ \ }^2$, respectively. The dual variables of $x$ and $t$ are
  denoted as $\xi$ and $\tau$, respectively. Thus $\widetilde{u}(\xi,\tau)
  =\mathcal{F}(u)(\xi,\tau)$, $\widehat{u}^1(\xi,t)=\mathcal{F}_1(u)(\xi,t)$,
  and $\widehat{u}^2(x,\tau)=\mathcal{F}_2(u)(x,\tau)$. In case no confusion
  may occur we often omit $1$ and $2$ in the notations $\widehat{\ \ }^1$ and
  $\widehat{\ \ }^2$, so that $\widehat{\varphi}(\xi)=\mathcal{F}_1(\varphi)
  (\xi)$ for $\varphi=\varphi(x)$, and $\widehat{\psi}(\tau)=\mathcal{F}_2
  (\psi)(\tau)$ for $\psi=\psi(t)$. The inverses of $\mathcal{F}$,
  $\mathcal{F}_1$ and $\mathcal{F}_2$ are denoted by $\mathcal{F}^{-1}$,
  $\mathcal{F}_1^{-1}$ and $\mathcal{F}_2^{-1}$, respectively.
\item $\mathcal{H}$ denotes the Hilbert transformation, i.e., $\mathcal{H}
  \varphi=\mathcal{F}_\xi^{-1}[-i\sgn\xi\cdot\widehat{\varphi}(\xi)]$ for
  $\varphi\in S'(\mathbb{R})$ such that $\sgn\xi\cdot\widehat{\varphi}(\xi)$
  makes sense and belongs to $S'(\mathbb{R})$. If $\varphi$ is a locally
  integrable function they we have
$$
  \mathcal{H}\varphi=
  \frac{1}{\pi}\mbox{p.v.}\!\int_{-\infty}^\infty\frac{\varphi(x)}{x-y}dy
$$
  in case the right-hand side makes sense.
\item For a real $s$, $D_x^s$ and $D_t^s$ denote absolute derivatives of
  order $s$ in the $x$ and $t$ variables, respectively, i.e., $D_x^s\varphi(x)
  =\mathcal{F}_1^{-1}[|\xi|^s\widehat{\varphi}(\xi)]$ for $\varphi\in
  S'(\mathbb{R})$ such that $|\xi|^s\widehat{\varphi}(\xi)\in S'(\mathbb{R})$,
  and similarly for $D_t^s$. $\langle D_x\rangle^s$ denotes the Fourier
  multiplier operator with symbol $\langle\xi\rangle^s=(1+|\xi|^2)^{s/2}$,
  i.e., $\langle D_x\rangle^s\varphi(x)=\mathcal{F}_1^{-1}[\langle\xi\rangle^s
  \widehat{\varphi}(\xi)]$ for $\varphi\in S'(\mathbb{R})$.
\item For a real $s$, $\dot{H}^{s}$ and $H^{s}$ respectively denote the
  homogeneous and inhomogeneous $L^2$-type Sobolev spaces on $\mathbb{R}$ of
  index $s$.
\end{itemize}

\section{Review of $L^2$ well-posedness}

  Before giving the proof of Theorem 1.1, let us first make a short review to
  the well-posedness result established in \cite{IK}. In some previous work
  (cf. \cite{I, P, T} for instance) it has been proved that the problem (1.1)
  is globally well-posed in $H^\sigma(\mathbb{R})$ for large $s$, and the best
  result is $\sigma\geq 1$ obtained by Tao in \cite{T}. By these results, there
  exists a continuous mapping $S^\infty:H^\infty(\mathbb{R}):=\cap_{\sigma\geq 0}
  H^\sigma(\mathbb{R})\to C(\mathbb{R},H^\infty(\mathbb{R}))$, such that for
  every $\phi\in H^\infty(\mathbb{R})$, $u=S^\infty\phi\in C^\infty(\mathbb{R},
  H^\infty(\mathbb{R}))$ is a solution of (1.1). For $T>0$ let $S_T^\infty:
  H^\infty(\mathbb{R})\to C([-T,T],H^\infty(\mathbb{R}))$ be the restriction of
  the mapping $S^\infty$ to the time interval $[-T,T]$. The result of \cite{IK}
  shows that the restriction $\sigma\geq 1$ can be weakened to $\sigma\geq 0$.
  We copy the main result of \cite{IK} (see Theorem 1.1 there) as follows:
\medskip

  {\bf Theorem 2.1}\ \ {\em $(a)$ For any $T>0$, the mapping $S_T^\infty:
  H^\infty(\mathbb{R})\to C([-T,T],H^\infty(\mathbb{R}))$ extends uniquely to
  a continuous mapping $S_T^0: L^2(\mathbb{R})\to C([-T,T],L^2(\mathbb{R}))$
  and $\|S_T^0(\phi)(\cdot,t)\|_2\\ =\|\phi\|_2$ for any $t\in [-T,T]$ and
  $\phi\in L^2(\mathbb{R})$. Moreover, for any $\phi\in L^2(\mathbb{R})$, the
  function $u=S_T^0(\phi)$ solves the initial-value problem $(1.1)$ in
  $C([-T,T],H^{-2}(\mathbb{R}))$.

  $(b)$ In addition, for any $\sigma\geq 0$, $S_T^0(H^\sigma(\mathbb{R}))
  \subseteq C([-T,T],H^\sigma(\mathbb{R}))$, $\|S_T^0(\phi)(\cdot,t)\|_{C([-T,T],
  H^\sigma)}\\ \leq C(T,\sigma,\|\phi\|_{H^\sigma})$, and the mapping $S_T^\sigma
  =S_T^0|_{H^\sigma(\mathbb{R})}: H^\sigma(\mathbb{R})\to C([-T,T],
  H^\sigma(\mathbb{R}))$ is continuous. $\quad\Box$}
\medskip

  {\bf Remark}\ \ From the discussion of \cite{IK} we see that for any $\phi\in
  L^2(\mathbb{R})$, the solution $u=S_T^0(\phi)$ has more regularity than merely
  being in $C([-T,T],L^2(\mathbb{R}))$; for instance, we have $u\in L^8([-T,T],
  L^4(\mathbb{R}))$ (cf. Lemma 3.6 in Section 3 below and note that
  $(4,8)$ is an admissible pair). Thus by inhomogeneous Strichartz estimates we
  see that $\int_0^t e^{-(t-t')\mathcal{H}\partial_x^2}u^2(\cdot,t') dt'\in
  C([-T,T],L^2(\mathbb{R}))\cap L^q([-T,T],L^p(\mathbb{R}))$ for any admissible
  pair $(p,q)$. Noticing this fact, it can be easily seen that for any $\phi\in
  L^2(\mathbb{R})$, $u=S_T^0(\phi)$ also solves the initial-value problem (1.1)
  in the sense that it satisfies the integral equation
$$
  u(\cdot,t)=e^{-t\mathcal{H}\partial_x^2}u_{0}+{1\over 2}\partial_x\!
  \int_0^t e^{-(t-t')\mathcal{H}\partial_x^2}u^2(\cdot,t') dt'
$$
  for $(x,t)\in\mathbb{R}\times (-T,T)$ in distribution sense. Conversely, it
  can also be easily seen that if a solution $u$ (in distribution sense) of
  this integral equation has certain regularity, for instance, $u\in C([-T,T],
  L^2(\mathbb{R}))\cap L^8([-T,T],L^4(\mathbb{R}))$, then $u$ also solves the
  initial-value problem (1.1) in $C([-T,T],H^{-2}(\mathbb{R}))$.
\medskip

  The main ingredients in proving the above result are a gauge transformation and
  the spaces $F^\sigma$ ($\sigma\geq 0$). For our purpose we review these
  ingredients in the following paragraphs.

  Let $P_{{\rm low}}$, $P_{\pm{\rm high}}$ and $P_{\pm}$ be projection operators
  on $L^2(\mathbb{R})$ defined respectively by
$$
  P_{{\rm low}}(\phi)=\mathcal{F}_1^{-1}(\widehat{\phi}
  \chi_{[-2^{10},2^{10}]}), \quad
  P_{\pm{\rm high}}(\phi)=\mathcal{F}_1^{-1}(\widehat{\phi}
  \chi_{\pm[2^{10},\infty)}),
$$
$$
  P_{\pm}(\phi)=\mathcal{F}_1^{-1}(\widehat{\phi}\chi_{\pm[0,\infty)}),
$$
  where $\chi_E$ (for given subset $E$ of $\mathbb{R}$) denotes the
  characteristic function of the subset $E$. Let $\phi\in H^\infty(\mathbb{R})$
  and set
$$
  \phi_{{\rm low}}=P_{{\rm low}}(\phi), \quad
  \phi_{\pm{\rm high}}=P_{\pm{\rm high}}(\phi).
$$
  It can be easily verified that for real-valued $\phi$, the function
  $\phi_{{\rm low}}$ is also real-valued. Let $u_0=S^\infty(\phi_{{\rm low}})$
  be the solution of the following problem:
$$
\left\{
\begin{array}{l}
  \partial_t u_0+\mathcal{H}\partial_x^2 u_0+
  \partial_x(u_0^2/2)=0, \quad
  x\in\mathbb{R},\;\;\; t\in\mathbb{R},\\
  u_0(x,0)=\phi_{{\rm low}}(x), \quad x\in\mathbb{R}.
\end{array}
\right.
\eqno{(2.1)}
$$
  Note that since $\phi_{{\rm low}}$ is real-valued, we have that $u_0$ is also
  real-valued. Besides, since $\|\phi_{{\rm low}}\|_{H^\sigma}\leq C_\sigma
  \|\phi\|_{L^2}$ for any $\sigma\geq 0$, it follows from the equation of $u_0$
  that
$$
  \sup_{|t|\leq T}\|\partial_t^{\sigma_1}\partial_x^{\sigma_2}
  u_0(\cdot,t)\|_{L_x^2}\leq
  C_{\sigma_1,\sigma_2}(\|\phi\|_{L^2})\|\phi\|_{L^2},
  \quad \sigma_1,\sigma_2\in \mathbb{Z}\cap [0,\infty).
\eqno{(2.2)}
$$
  We define a gauge $U_0$ as follows: First let $U_0(0,t)$ be the solution of
  the following problem:
$$
  \partial_t U_0(0,t)+{1\over 2}\mathcal{H}\partial_x u_0(0,t)+
  {1\over 4}u_0^2(0,t)=0 \quad \mbox{for}\;\; t\in\mathbb{R},
  \quad \mbox{and}\;\; U_0(0,0)=0,
$$
  and next extend $U_0(x,t)$ to all $x\in\mathbb{R}$ (for fixed $t\in
  \mathbb{R})$) by using the following equation:
$$
  \partial_x U_0(x,t)={1\over 2}u_0(x,t).
$$
  Note that since $u_0$ is real-valued, we see that $U_0$ is also real-valued.
  Besides, for any integers $\sigma_1,\sigma_2\geq 0$, $(\sigma_1,\sigma_2)\neq
  (0,0)$,
$$
  \sup_{|t|\leq T}\|\partial_t^{\sigma_1}\partial_x^{\sigma_2}
  U_0(\cdot,t)\|_{L_x^2}\leq
  C_{\sigma_1,\sigma_2}(\|\phi\|_{L^2})\|\phi\|_{L^2}.
\eqno{(2.3)}
$$
  We now define
$$
\left\{
\begin{array}{l}
  w_+=e^{iU_0}P_{+{\rm high}}(u-u_0),\\
  w_-=e^{-iU_0}P_{-{\rm high}}(u-u_0),\\
  w_0=P_{{\rm low}}(u-u_0).
\end{array}
\right.
$$
  Then $(w_+,w_-,w_0)$ satisfies the following system of equations (see
  (2.10), (2.12) and (2.14) of \cite{IK}):
$$
\left\{
\begin{array}{l}
  \partial_tw_++\mathcal{H}\partial_x^2w_+=E_+(w_+,w_-,w_0), \quad
  x\in\mathbb{R},\;\;\; t\in\mathbb{R},\\
  \partial_tw_-+\mathcal{H}\partial_x^2w_-=E_-(w_+,w_-,w_0), \quad
  x\in\mathbb{R},\;\;\; t\in\mathbb{R},\\
  \partial_tw_0+\mathcal{H}\partial_x^2w_0=E_0(w_+,w_-,w_0), \quad
  x\in\mathbb{R},\;\;\; t\in\mathbb{R},\\
  (w_+,w_-,w_0)|_{t=0}=(e^{iU_0(\cdot,0)}\phi_{+{\rm high}},
  e^{-iU_0(\cdot,0)}\phi_{-{\rm high}},0),
\end{array}
\right.
\eqno{(2.4)}
$$
  where (see (2.11), (2.13) and (2.15) of \cite{IK})
$$
\begin{array}{rl}
  E_+(w_+,&w_-,w_0)=-e^{iU_0}P_{+{\rm high}}[\partial_x(e^{-iU_0}w_+
  +e^{iU_0}w_-+w_0)^2/2]
\\ [0.2cm]
  &-e^{iU_0}P_{+{\rm high}}\{\partial_x[u_0\cdot
  P_{-{\rm high}}(e^{iU_0}w_-)+u_0\cdot P_{{\rm low}}(w_0)]\}
\\ [0.2cm]
  &+e^{iU_0}(P_{-{\rm high}}+P_{{\rm low}})\{\partial_x[u_0\cdot
  P_{+{\rm high}}(e^{-iU_0}w_+)]\}
\\ [0.2cm]
  &+2iP_{-}\{\partial_x^2[e^{iU_0}P_{+{\rm high}}(e^{-iU_0}w_+)]\}
\\ [0.2cm]
  &-P_{+}(\partial_x u_0)\cdot w_+,
\end{array}
$$
$$
\begin{array}{rl}
  E_-(w_+,&w_-,w_0)=-e^{-iU_0}P_{-{\rm high}}[\partial_x(e^{-iU_0}w_+
  +e^{iU_0}w_-+w_0)^2/2]
\\ [0.2cm]
  &-e^{-iU_0}P_{-{\rm high}}\{\partial_x[u_0\cdot
  P_{+{\rm high}}(e^{-iU_0}w_+)+u_0\cdot P_{{\rm low}}(w_0)]\}
\\ [0.2cm]
  &+e^{-iU_0}(P_{+{\rm high}}+P_{{\rm low}})\{\partial_x[u_0\cdot
  P_{-{\rm high}}(e^{iU_0}w_-)]\}
\\ [0.2cm]
  &-2iP_{+}\{\partial_x^2[e^{-iU_0}P_{-{\rm high}}(e^{iU_0}w_-)]\}
\\ [0.2cm]
  &-P_{-}(\partial_x u_0)\cdot w_-,
\end{array}
$$
$$
  E_0(w_+,w_-,w_0)=-\frac{1}{2}P_{{\rm low}}\{\partial_x[(e^{-iU_0}w_+
  +e^{iU_0}w_-+w_0+u_0)^2-u_0^2]\}.
$$
  It is immediate to see that the following relation holds (see Lemma 2.1
  of \cite{IK}):
$$
  u=e^{-iU_0}w_+ + e^{iU_0}w_- + w_0 + u_0.
\eqno{(2.5)}
$$

  The mapping $u\to (w_+,w_-,w_0)$ is called {\em gauge transform} (in more
  precise sense the components $u_0$ and $U_0$ should also be comprised into
  this notion; but for simplicity of the notation we omit them).
  The above deduction shows that if $u$ is a smooth solution of (1.1) (or
  more precisely, a solution of (1.1) whose initial data belong to
  $H^\infty(\mathbb{R})$) then $(w_+,w_-,w_0)$ is a solution of (2.4). The
  converse assertion cannot be directly verified. The proof (for smooth $\phi$)
  that if $(w_+,w_-,w_0)$ is a solution of (2.4) then the expression $u$ given
  by (2.5) is a solution of (1.1) is given in Section 10 of \cite{IK};
  see (10.38) in \cite{IK}. The main idea in the proof of Theorem 3.1 is as
  follows: First one proves that for small initial data the problem (2.4) is
  well-posed in suitable function spaces; in particular it has a solution in
  $C([-T,T],(L^2(\mathbb{R}))^3)$ depending continuously on the initial data.
  Using this fact and the relation (2.5) established for smooth solutions, one
  then proves that the solution operator $S_T^\infty$ defined for smooth data
  can be extended into a continuous mapping from $L^2(\mathbb{R})$ to $C([-T,T],
  L^2(\mathbb{R}))$. Since $H^\infty(\mathbb{R})$ is dense in $L^2(\mathbb{R})$,
  the extension is unique, and is denoted by $S_T^0$. For any given $T>0$ and
  $\phi\in L^2(\mathbb{R})$ with sufficiently small norm, $u=S_T^0(\phi)$ then
  defines a solution of the problem (1.1) for $|t|\leq T$. A standard scaling
  argument enables one to convert this small-data assertion into a local
  well-posedness result for (1.1) for arbitrary initial data in
  $L^2(\mathbb{R})$, and the $L^2$ conservation law then yields the desired
  global well-posedness result.

  Well-posedness of the problem (2.4) is established in a class of spaces
  $F^\sigma$ ($\sigma\geq 0$), whose definition is given below. Let $\eta_0:
  \mathbb{R}\to [0,1]$ denote an even function supported in $[-8/5,8/5]$
  and equal to $1$ in $[-5/4,5/4]$. For $k\in\mathbb{Z}$, $k\geq 1$, let
  $\eta_k(\xi)=\eta_0(2^{-k}\xi)-\eta_0(2^{-k+1}\xi)$. We also denote, for
  all $k\in\mathbb{Z}$, $\chi_k(\xi)=\eta_0(2^{-k}\xi)-\eta_0(2^{-k+1}\xi)$.
  It follows that
$$
  \sum_{k=0}^\infty\eta_k(\xi)=1\quad \mbox{for}\;\;
  \xi\in\mathbb{R},
$$
  and
$$
  \sum_{k=-\infty}^\infty\chi_k(\xi)=1\quad \mbox{for}\;\;
  \xi\in\mathbb{R}\backslash\{0\}.
$$
  Note that $\supp\chi_k\subseteq [-(8/5)2^k,-(5/8)2^k]\cup [(5/8)2^k,
  (8/5)2^k]$ for all $k\in\mathbb{Z}$, and $\supp\eta_k\subseteq [-(8/5)2^k,
  -(5/8)2^k]\cup [(5/8)2^k,(8/5)2^k]$ for $k\geq 1$. For $k\in\mathbb{Z}$ we
  denote $I_k=[-2^{k+1},-2^{k-1}]\cup [2^{k-1},2^{k+1}]$, and for $k\in
  \mathbb{Z}$, $k\geq 1$, we also denote $\tilde{I}_k=[-2,2]$ if $k=0$ and
  $\tilde{I}_k=I_k$ if $k\geq 1$. Next, we denote
$$
  \omega(\xi)=-\xi|\xi| \;\; (\xi\in\mathbb{R}) \quad \mbox{and} \quad
  \beta_{k,j}=1+2^{(j-2k)/2} \;\;(j,k\in\mathbb{Z}),
$$
  and for $k\in\mathbb{Z}$ and $j\geq 0$ we let
$$
  D_{k,j}=\left\{
\begin{array}{l}
  \{(\xi,\tau)\in\mathbb{R}\times\mathbb{R}:\xi\in I_k,\,
    \tau\!-\!\omega(\xi)\in\tilde{I}_j\}\;\;\; \mbox{if}\;\; k\geq 1;
\\ [0.2cm]
  \{(\xi,\tau)\in\mathbb{R}\times\mathbb{R}:\xi\in I_k,\,
    \tau\in\tilde{I}_j\}\;\;\; \mbox{if}\;\; k\leq 0.
\end{array}
\right.
$$
  We now define spaces $\{Z_k\}_{k=0}^\infty$ as follows:
$$
  Z_k=\left\{
\begin{array}{ll}
    X_k+Y_k \quad & \mbox{if}\;\; k=0\;\;\mbox{or}\;\; k\geq 100,
\\ [0.2cm]
     X_k \quad & \mbox{if}\;\; 1\leq k\leq 99,
\end{array}
\right.
$$
  where
$$
\begin{array}{rl}
  X_k=&\{f\in L^2(\mathbb{R}\times\mathbb{R}): f\;\;\mbox{supported in}\;\;
  I_k\times\mathbb{R}\;\; \mbox{and}\\ [0.2cm]
  &\|f\|_{X_k}:=\sum_{j=0}^\infty 2^{j/2}\beta_{k,j}\|
   \eta_j(\tau\!-\!\omega(\xi))f(\xi,\tau)\|_{L^2_{\xi,\tau}}<\infty\}
   \quad \mbox{for}\;\; k\geq 1,
\\ [0.2cm]
  X_0=&\{f\in L^2(\mathbb{R}\times\mathbb{R}): f\;\;\mbox{supported in}\;\;
  \tilde{I}_0\times\mathbb{R}\;\; \mbox{and}\\ [0.2cm]
  &\|f\|_{X_0}:=\sum_{j=0}^\infty\sum_{l=-\infty}^1 2^{j-l}\|\eta_j(\tau)
   \chi_l(\xi)f(\xi,\tau)\|_{L^2_{\xi,\tau}}<\infty\},
\end{array}
$$
  and
$$
\begin{array}{rl}
  Y_k=&\{f\in L^2(\mathbb{R}\times\mathbb{R}): f\;\;\mbox{supported in}\;\;
    \cup_{j=0}^{k-1}D_{k,j}\;\; \mbox{and}\\ [0.2cm]
  &\|f\|_{Y_k}:=2^{-k/2}\|\mathcal{F}^{-1}[(\tau\!-\!\omega(\xi)\!+\!i)
   f(\xi,\tau)]\|_{L_x^1L_t^2}
   <\infty\} \quad \mbox{for}\;\; k\geq 1,
\\ [0.2cm]
  Y_0=&\{f\in L^2(\mathbb{R}\times\mathbb{R}): f\;\;\mbox{supported in}\;\;
  \tilde{I}_0\times\mathbb{R}\;\; \mbox{and}\\ [0.2cm]
  &\|f\|_{X_0}:=\sum_{j=0}^\infty 2^{j}\|\mathcal{F}^{-1}[\eta_j(\tau)
   f(\xi,\tau)]\|_{L_x^1L_t^2}<\infty\}.
\end{array}
$$
  Let $\sigma\geq 0$. The space $F^\sigma$ is defined as follows:
$$
  F^\sigma=\{u\in S'(\mathbb{R}\times\mathbb{R}):\|u\|_{F^\sigma}^2
    :=\sum_{k=0}^\infty 2^{2\sigma k}\|\eta_k(\xi)(I-\partial_\tau^2)
    \widetilde{u}(\xi,\tau)\|_{Z_k}^2<\infty\}.
$$
  $F^\sigma$ is the space which plays a role in the study of well-posedness of
  the problem (2.4) similar to the role of the Bourgain space $X^{\sigma,b}$
  in the study of well-posedness of the KdV equation. However, the
  corresponding space in the space variable is not $H^\sigma$, but instead
  $\widetilde{H}^\sigma$, which is defined as follows. First we define
$$
\begin{array}{rl}
  B_0=&\{f\in L^2(\mathbb{R}): f\;\;\mbox{supported in}\;\;
  \tilde{I}_0\;\; \mbox{and}\\
  &\displaystyle\|f\|_{B_0}:=\inf_{f=g+h}\Big\{\|F_1^{-1}(g)\|_{L_x^1}
   +\sum_{l=-\infty}^1 2^{-l}\|\chi_l h\|_{L_\xi^2}\Big\}<\infty\}.
\end{array}
$$
  It is clear that $\|f\|_{L^2}\leq 2\|f\|_{B_0}$. Then we define
$$
  \widetilde{H}^\sigma=\{\phi\in L^2(\mathbb{R}):
  \|\phi\|_{\widetilde{H}^\sigma}^2:=\|\eta_0\widehat{\phi}\|_{B_0}^2+
  \sum_{k=1}^\infty 2^{2\sigma k}\|\eta_k\widehat{\phi}\|_{L^2}^2<\infty\}.
$$
  Since $B_0\hookrightarrow L^2$, we see that $\widetilde{H}^\sigma
  \hookrightarrow {H}^\sigma$, and $\|\phi\|_{{H}^\sigma}\leq\sqrt{2}
  \|\phi\|_{\widetilde{H}^\sigma}$. By Lemma 4.2 of \cite{IK} we know that
$$
  F^\sigma\subseteq C(\mathbb{R},\widetilde{H}^\sigma) \quad
  \mbox{for any}\;\;\sigma\geq 0,
$$
  and the embedding is continuous.

  Given $T>0$, we denote by $F^\sigma_T$ the restriction of the space
  $F^\sigma$ in $\mathbb{R}\times [-T,T]$. From the discussion in Section 10
  of \cite{IK} we have:
\medskip

  {\bf Theorem 2.2}\ \ {\it Given $T>0$, there exists corresponding $\epsln>0$
  such that for any $\psi_+$, $\psi_-$, $\psi_0\in\widetilde{H}^0$ satisfying
$$
  \|\psi_+\|_{\widetilde{H}^0}+\|\psi_-\|_{\widetilde{H}^0}+
  \|\psi_0\|_{\widetilde{H}^0}\leq\epsln,
\eqno{(2.6)}
$$
  the initial value problem
$$
\left\{
\begin{array}{l}
  \partial_tv_++\mathcal{H}\partial_x^2v_+=E_+(v_+,v_-,v_0), \quad
  x\in\mathbb{R},\;\;\; t\in [-T,T],\\
  \partial_tv_-+\mathcal{H}\partial_x^2v_-=E_-(v_+,v_-,v_0), \quad
  x\in\mathbb{R},\;\;\; t\in [-T,T],\\
  \partial_tv_0+\mathcal{H}\partial_x^2v_0=E_0(v_+,v_-,v_0), \quad
  x\in\mathbb{R},\;\;\; t\in [-T,T],\\
  (v_+,v_-,v_0)|_{t=0}=(\psi_+,\psi_-,\psi_0)
\end{array}
\right.
\eqno{(2.7)}
$$
  has a solution $(v_+,v_-,v_0)\in (F^0_T)^3$, which lies in a small ball
  $B_{\epsln'}$ of $(F^0_T)^3$ and is the unique solution of $(2.7)$ in this
  ball, where $\epsln'=\epsln'(\epsln)>0$ is such that $\epsln'\to 0$ as
  $\epsln\to 0$, and the mapping
  $(\psi_+,\psi_-,\psi_0)\to (v_+,v_-,v_0)$ from $(\widetilde{H}^0)^3$ to
  $(F^0_T)^3$ is continuous. Moreover, if $\psi_+$, $\psi_-$, $\psi_0\in
  \widetilde{H}^\sigma$ for some $\sigma>0$ then $(v_+,v_-,v_0)\in
  (F^\sigma_T)^3$, and the mapping $(\psi_+,\psi_-,\psi_0)\to (v_+,v_-,v_0)$
  from $(\widetilde{H}^\sigma)^3$ to $(F^\sigma_T)^3$ is continuous.
  $\quad\Box$}
\medskip

  As we saw before, the relation (2.5) connecting the solution $u$ of (1.1)
  with the solution $(w_+,w_-,w_0)$ of (2.4) was only established for smooth
  initial data. With the aid of Theorem 2.2, we can extend it to all solutions
  with $L^2$ data, i.e., we have the following result:
\medskip

  {\bf Theorem 2.3}\ \ {\it Given $T>0$, there exists corresponding $\epsln>0$
  such that for any $\phi\in L^2(\mathbb{R})$ satisfying $\|\phi\|_{L^2}\leq
  \epsln$, the solution $u=S_T^0(\phi)$ of the problem $(1.1)$ has the
  expression $(2.5)$, with $u_0$ and $U_0$ as in $(2.1)$--$(2.3)$, and $(w_+,
  w_-,w_0)$ being the unique solution of $(2.4)$ in $($a small neighborhood of
  the origin of$)$ the space $(F^0_T)^3$ with norm $\leq\epsln'(\epsln)$.}
\medskip

  {\it Proof}:\ \ By (2.9) and Lemma 10.1 of \cite{IK} we see that for any
  $\phi\in H^\sigma$ ($\sigma\geq 0$) we have
$$
  (e^{iU_0(\cdot,0)}\phi_{+{\rm high}},e^{-iU_0(\cdot,0)}
  \phi_{-{\rm high}},0)\in (\widetilde{H}^\sigma)^3,
$$
  and the mapping $\phi\to (\psi_+(\phi),\psi_-(\phi),\psi_0(\phi)):=
  (e^{iU_0(\cdot,0)}\phi_{+{\rm high}},e^{-iU_0(\cdot,0)}\phi_{-{\rm high}},0)$
  from $H^\sigma$ to $(\widetilde{H}^\sigma)^3$ is continuous. Using this
  assertion particularly to $\sigma=0$, we see that for $\epsln>0$ as in (2.6),
  there exists corresponding $\epsln'>0$ such that if $\|\phi\|_{L^2}\leq
  \epsln'$ then
$$
  \|\psi_+(\phi)\|_{\widetilde{H}^0}+\|\psi_-(\phi)\|_{\widetilde{H}^0}+
  \|\psi_0(\phi)\|_{\widetilde{H}^0}\leq\epsln.
$$
  By Theorem 2.2, for such $\phi\in L^2(\mathbb{R})$ the problem (2.4) has a
  unique solution $(w_+,w_-,w_0)\in (F^0_T)^3$. We now assume that $\phi\in
  L^2(\mathbb{R})$ is a such function, i.e., $\|\phi\|_{L^2}\leq\epsln'$, and
  let $u=S_T^0(\phi)$. Let $\phi_n=\mathcal{F}_1^{-1}(\widehat{\phi}
  \chi_{[-2^{10}n,2^{10}n]})$, $n=1,2,\cdots$. Then we have $\phi_n\in
  H^\infty(\mathbb{R})$,
$$
  \|\phi_n\|_{L^2}\leq\|\phi\|_{L^2}\leq\epsln', \;\;\; n=1,2,\cdots, \quad
  \mbox{and} \quad \lim_{n\to\infty}\|\phi_n-\phi\|_{L^2}=0.
$$
  Let $u_n=S_T^\infty(\phi_n)$, and let $u_{n0}$, $U_{n0}$, $w_{n+}$, $w_{n-}$,
  $w_{n0}$ be the corresponding counterparts of $u_{0}$, $U_{0}$, $w_{+}$,
  $w_{-}$, $w_{0}$ defined before when $\phi$ is replaced by $\phi_n$, $n=1,2,
  \cdots$. Then we have
$$
  u_n=e^{-iU_{n0}}w_{n+}+e^{iU_{n0}}w_{n-}+w_{n0}+u_{n0},
  \quad n=1,2,\cdots.
$$
  From the special construction of the function $\phi_n$ we see that
  $P_{{\rm low}}(\phi_n)=P_{{\rm low}}(\phi)$ for all $n\in\mathbb{N}$, so that
  $u_{n0}=u_0$ for all $n\in\mathbb{N}$ and, consequently, $U_{n0}=U_0$ for all
  $n\in\mathbb{N}$. Thus, the above relations can be rewritten as follows:
$$
  u_n=e^{-iU_0}w_{n+}+e^{iU_0}w_{n-}+w_{n0}+u_0,
  \quad n=1,2,\cdots.
\eqno{(2.8)}
$$
  Note that $(w_{n+},w_{n-},w_{n0})$ are in a small ball in $(F_T^0)^3$.
  Using Lemma 10.1 in \cite{IK} and the facts that $U_{n0}=U_0$ for all $n\in
  \mathbb{N}$ and $\phi_n\to\phi$ strongly in $L^2(\mathbb{R})$, we see that
$$
  \|\psi_+(\phi_n)-\psi_+(\phi)\|_{\widetilde{H}^0}
  +\|\psi_-(\phi_n)-\psi_-(\phi)\|_{\widetilde{H}^0}
  +\|\psi_0(\phi_n)-\psi_0(\phi)\|_{\widetilde{H}^0}\to 0
  \quad \mbox{as}\;\; n\to\infty.
$$
  Thus, by the continuity assertion in Theorem 2.2 we conclude that
$$
  \|w_{n+}-w_+\|_{F_T^0}+\|w_{n-}-w_-\|_{F_T^0}+
  \|w_{n0}-w_0\|_{F_T^0}\to 0 \quad \mbox{as}\;\; n\to\infty.
$$
  Since $F_T^0$ is continuously embedded into $C([-T,T],L^2(\mathbb{R}))$,
  this implies that
$$
  \sup_{|t|\leq T}\|w_{n+}(\cdot,t)-w_+(\cdot,t)\|_2
  +\sup_{|t|\leq T}\|w_{n-}(\cdot,t)-w_-(\cdot,t)\|_2
  +\sup_{|t|\leq T}\|w_{n0}(\cdot,t)-w_0(\cdot,t)\|_2\to 0
  \quad \mbox{as}\;\; n\to\infty.
$$
  Hence, by letting $n\to\infty$ in (2.8) and using the facts that $u=
  S_T^0(\phi)=\lim_{n\to\infty}S_T^\infty(\phi_n)$ (in $C([-T,T],
  L^2(\mathbb{R}))$ norm) and $u_n=S_T^\infty(\phi_n)$, we see that (2.5)
  follows. To get the desired assertion we only need to re-denote $\epsln'$
  as $\epsln$. This completes the proof. $\quad\Box$

\section{Preliminary estimates}

\hskip 2em
  {\bf  Lemma 3.1}\ \ {\it For any $k\geq 0$, if $f_k\in Z_k$ then}
$$
  \|\mathcal{F}^{-1}(f_k)\|_{L_x^\infty L_t^2}\leq C 2^{-k/2}\|f_k\|_{Z_k}.
\eqno{(3.1)}
$$

  {\it Proof}:\ \  For $k\geq 1$, this assertion has been proved in \cite{IK}
 (see Lemma 4.2 $(c)$ of \cite{IK}). Hence, in the sequel we only consider the
  case $k=0$.

  Let $\phi_0\in C^\infty_0(\mathbb{R})$ such that $\phi_0(\xi)=1$ for $|\xi|
  \leq 2$. Next let $\psi\in C^\infty_0(\mathbb{R})$ such that $\supp\psi
  \subseteq [-5/2,5/2]$, $\psi(\xi)=1$ for $|\xi|\leq 2/5$, and define
  $\phi_k(\tau)=\psi(2^{-k-2}\tau)-\psi(2^{-k+2}\tau)=\psi_0(2^{-k}
  \tau)$ for $k\geq 1$, where $\psi_0(\tau)=\psi(\tau/4)-\psi(4\tau)$.
  Then $\phi_k\in C^\infty_0(\mathbb{R})$ and $\phi_k(\tau)=1$ if $(5/8)2^k
  \leq|\tau|\leq (8/5)2^k$ ($k\geq 1$). Hence, since $\supp\eta_k\subseteq
  [-(8/5)2^k,-(5/8)2^k]\cup [(5/8)2^k,(8/5)2^k]$ for $k\geq 1$, we have
  $\phi_k(\tau)\eta_k(\tau)=\eta_k(\tau)$ for all $k\geq 1$.

  We first assume that $f_0\in X_0$. Then, since $f_0$ is supported in
  $\widetilde{I}\times\mathbb{R}$, we have
$$
\begin{array}{rl}
  f_0(\xi,\tau)=&\displaystyle\phi_0(\xi)f_0(\xi,\tau)=\sum_{j=0}^\infty
  \sum_{l=-\infty}^1\eta_j(\tau)\chi_l(\xi)\phi_0(\xi)f_0(\xi,\tau)
\\[0.3cm]
  =&\displaystyle\sum_{j=0}^\infty\sum_{l=-\infty}^1\phi_j(\tau)
  \eta_j(\tau)\chi_l(\xi)\phi_0(\xi)f_0(\xi,\tau)
 =\sum_{j=0}^\infty\sum_{l=-\infty}^1\phi_0(\xi)\phi_j(\tau)\cdot
  f_{0jl}(\xi,\tau),
\end{array}
$$
  where $f_{0jl}(\xi,\tau)=\eta_j(\tau)\chi_l(\xi)f_0(\xi,\tau)$. Hence,
\setcounter{equation}{1}
\begin{eqnarray}
  \mathcal{F}^{-1}[f_0(\xi,\tau)]&=&\displaystyle\sum_{j=0}^\infty
  \sum_{l=-\infty}^1
  \mathcal{F}^{-1}[\phi_0(\xi)\phi_j(\tau)\cdot f_{0jl}(\xi,\tau)]
\nonumber \\
  &=&\displaystyle\sum_{j=0}^\infty\sum_{l=-\infty}^1
  \mathcal{F}^{-1}[\phi_0(\xi)\phi_j(\tau)]\ast
  \mathcal{F}^{-1}[f_{0jl}(\xi,\tau)],
\end{eqnarray}
which yields
$$
\begin{array}{rl}
  \|\mathcal{F}^{-1}(f_0)\|_{L_x^\infty L_t^2}\leq
  &\displaystyle\sum_{j=0}^\infty\sum_{l=-\infty}^1
  \|\mathcal{F}^{-1}[\phi_0(\xi)\phi_j(\tau)]\|_{L_x^2 L_t^1}
  \|\mathcal{F}^{-1}[f_{0jl}(\xi,\tau)]\|_{L_{x,t}^2}
\\ [0.3cm]
  \leq & C\displaystyle\sum_{j=0}^\infty\sum_{l=-\infty}^1
  \|\phi_0(\xi)\|_{L_\xi^2}\|2^j\tilde{\psi_0}(2^jt)]\|_{L_t^1}
  \|f_{0jl}(\xi,\tau)\|_{L_{\xi,\tau}^2}
 \qquad (\tilde{\psi_0}=\mathcal{F}_2^{-1}(\psi_0))
\\ [0.3cm]
  \leq & C\displaystyle\sum_{j=0}^\infty\sum_{l=-\infty}^1
  \|\eta_j(\tau)\chi_l(\xi)f_0(\xi,\tau)\|_{L_{\xi,\tau}^2}
  \leq C\|f_0\|_{X_0}.
\end{array}
$$
  We next assume that $f_0\in Y_0$. Then as before we have
$$
\begin{array}{rl}
  f_0(\xi,\tau)=&\displaystyle\phi_0(\xi)f_0(\xi,\tau)=\sum_{j=0}^\infty
  \eta_j(\tau)\phi_0(\xi)f_0(\xi,\tau)
  =\sum_{j=0}^\infty\phi_0(\xi)\phi_j(\tau)\cdot
  f_{0j}(\xi,\tau),
\end{array}
$$
  where $f_{0j}(\xi,\tau)=\eta_j(\tau)f_0(\xi,\tau)$. Hence,
\setcounter{equation}{2}
\begin{eqnarray}
  \mathcal{F}^{-1}[f_0(\xi,\tau)]&=&\displaystyle\sum_{j=0}^\infty
  \mathcal{F}^{-1}[\phi_0(\xi)\phi_j(\tau)\cdot f_{0j}(\xi,\tau)]
  =\sum_{j=0}^\infty
  \mathcal{F}^{-1}[\phi_0(\xi)\phi_j(\tau)]\ast
  \mathcal{F}^{-1}[f_{0j}(\xi,\tau)], \qquad
\end{eqnarray}
  which yields
$$
\begin{array}{rl}
  \|\mathcal{F}^{-1}(f_0)\|_{L_x^\infty L_t^2}\leq
  &\displaystyle\sum_{j=0}^\infty
  \|\mathcal{F}^{-1}[\phi_0(\xi)\phi_j(\tau)]\|_{L_x^\infty L_t^1}
  \|\mathcal{F}^{-1}[f_{0j}(\xi,\tau)]\|_{L_x^1L_t^2}
\\ [0.3cm]
  \leq & C\displaystyle\sum_{j=0}^\infty
\|\phi_0(\xi)\|_{L_\xi^1}\|2^j\tilde{\psi_0}(2^jt)]\|_{L_t^1}
  \|\mathcal{F}^{-1}[\eta_j(\tau)f_0(\xi,\tau)]\|_{L_x^1L_t^2}
 \qquad (\tilde{\psi_0}=\mathcal{F}_2^{-1}(\psi_0))
\\ [0.3cm]
  \leq & C\displaystyle\sum_{j=0}^\infty
  \|\mathcal{F}^{-1}[\eta_j(\tau)f_0(\xi,\tau)]\|_{L_x^1L_t^2}
  \leq C\|f_0\|_{Y_0}.
\end{array}
$$
  Now let $f_0\in Z_0$. Then there exists $g_0\in X_0$ and $h_0\in Y_0$
  such that
$$
  f_0=g_0+h_0 \quad \mbox{and} \quad
  \|g_0\|_{X_0}+\|h_0\|_{Y_0}\leq 2\|f_0\|_{Z_0}.
$$
  Thus
$$
  \|\mathcal{F}^{-1}(f_0)\|_{L_x^\infty L_t^2}
  \leq \|\mathcal{F}^{-1}(g_0)\|_{L_x^\infty L_t^2}\!+\!
  \|\mathcal{F}^{-1}(h_0)\|_{L_x^\infty L_t^2}
  \leq  C(\|g_0\|_{X_0}\!+\!\|h_0\|_{Y_0})
  \leq  C\|f_0\|_{Z_0}.
$$
  This completes the proof of  Lemma 3.1. $\quad\Box$
\medskip

  In the proof of the following lemma we shall use the following fact:
  If $f_k\in Z_k$ ($k\geq 0$) then for any $\alpha\geq 0$,
$$
  \||\xi|^\alpha f_k(\xi,\tau)\|_{Z_k}\leq C2^{\alpha k}
  \|f_k\|_{Z_k}.
$$
  This is an immediate consequence of Lemma 4.1 $a)$ of \cite{IK}.
\medskip

  {\bf  Lemma 3.2}\ \ {\it If $w\in F^0$ then for any $0\leq\theta<1$
  we have $D_x^{{\theta\over 2}}w\in L_x^\infty L_t^2$,}
  and
$$
  \|D_x^{{\theta\over 2}}w\|_{L_x^\infty L_t^2}\leq C\|w\|_{F^0}.
\eqno{(3.4)}
$$

  {\it Proof}:\ \ Let $f_k(\xi,\tau)=\eta_k(\xi)\widetilde{w}(\xi,\tau)$,
  $k=0,1,2,\cdots$, where $\widetilde{w}=F(w)$. Then $w\in F^0$ implies
  that $(I-\partial_\tau^2)f_k\in Z_k$, $k=0,1,2,\cdots$, and
$$
  \|w\|_{F^0}=\Big(\sum_{k=0}^\infty\|(I-\partial_\tau^2)f_k\|_{Z_k}^2
  \Big)^{{1\over 2}}<\infty.
$$
  Since $\widetilde{w}(\xi,\tau)=\displaystyle\sum_{k=0}^\infty
  f_k(\xi,\tau)$, for any $0\leq\theta<1$ we have
$$
  (1+t^2)D_x^{{\theta\over 2}}w(x,t)=\sum_{k=0}^\infty
  \mathcal{F}^{-1}[|\xi|^{{\theta\over 2}}(I-\partial_\tau^2)f_k(\xi,\tau)].
$$
  Hence, by  Lemma 3.1 we have
$$
\begin{array}{rl}
  &\|(1+t^2)D_x^{{\theta\over 2}}w(x,t)\|_{L_x^\infty L_t^2}\leq
  \displaystyle\sum_{k=0}^\infty\|\mathcal{F}^{-1}[|\xi|^{{\theta\over 2}}
  (I-\partial_\tau^2)f_k(\xi,\tau)]\|_{L_x^\infty L_t^2}
\\ [0.3cm]
  \leq & C\displaystyle\sum_{k=0}^\infty 2^{-{k\over 2}}
  \||\xi|^{{\theta\over 2}}(I-\partial_\tau^2)f_k(\xi,\tau)\|_{Z_k}
  \leq C\displaystyle\sum_{k=0}^\infty 2^{-{k\over 2}}
  2^{{\theta k\over 2}}\|(I-\partial_\tau^2)f_k(\xi,\tau)\|_{Z_k}
\\ [0.3cm]
  \leq & C\displaystyle\Big(\sum_{k=0}^\infty 2^{-(1\!-\!\theta)k}
  \Big)^{{1\over 2}}\Big(\sum_{k=0}^\infty\|(I-\partial_\tau^2)
  f_k(\xi,\tau)\|_{Z_k}^2\Big)^{{1\over 2}}
  \leq C\|w\|_{F^0}.
\end{array}
$$
  From this estimate (3.4) follows immediately. $\quad\Box$
\medskip

  {\bf Lemma 3.3}\ \ {\it Let $0\leq\theta<1$. For any $k\geq 0$, if $f_k\in
  Z_k$ then}
$$
  \|D_t^{{\theta\over 4}}\mathcal{F}^{-1}(f_k)\|_{L_x^\infty L_t^2}\leq
  C2^{-(1-\theta)k/2}\|f_k\|_{Z_k}.
\eqno{(3.5)}
$$

  {\it Proof}:\ \ We first assume that $k\geq 1$ and $f_k\in X_k$. Let
  $f_{k,j}(\xi,\tau)=\eta_j(\tau-\omega(\xi))f_k(\xi,\tau)$, $j\in
  \mathbb{Z}$, $j\geq 0$. Then $\supp f_{k,j}\subseteq D_{k,j}$,
$$
  f_k(\xi,\tau)=\sum_{j=0}^\infty f_{k,j}(\xi,\tau),
\eqno{(3.6)}
$$
  and
$$
  \|f_k\|_{X_k}=\sum_{j=0}^\infty 2^{j/2}(1+2^{(j-2k)/2})
  \|f_{k,j}(\xi,\tau)\|_{L^2_{\xi,\tau}}.
\eqno{(3.7)}
$$
  From (3.6) we have
$$
  D_t^{{\theta\over 4}}\mathcal{F}^{-1}(f_k)=
  \mathcal{F}^{-1}[|\tau|^{{\theta\over 4}}f_k(\xi,\tau)]=\sum_{j=0}^\infty
  \mathcal{F}^{-1}[|\tau|^{{\theta\over 4}}f_{k,j}(\xi,\tau)]
$$
  From the proof of (3.1) (see Line 3, Page 763 of \cite{IK}) we know
  that
$$
  \|\mathcal{F}^{-1}[|\tau|^{{\theta\over 4}}f_{k,j}(\xi,\tau)]\|_{L_x^\infty L_t^2}
  \leq C2^{-k/2}2^{j/2}\||\tau|^{{\theta\over 4}}
  f_{k,j}(\xi,\tau)\|_{L_{\xi,\tau}^2}.
$$
  Since $f_{k,j}$ is supported in $D_{k,j}$, we have $|\xi|\leq C2^k$ and
  $|\tau-\omega(\xi)||\leq C2^j$ for $(\xi,\tau)\in\supp(f_{k,j})$. If $j\leq
  2k$ then we have
$$
  |\tau|\leq |\tau-\omega(\xi)|+|\omega(\xi)|\leq
   C2^j+ C2^{2k}\leq C2^{2k},
$$
  so that
\begin{eqnarray}
  &&\|\mathcal{F}^{-1}[|\tau|^{{\theta\over 4}}f_{k,j}(\xi,\tau)]\|_{L_x^\infty L_t^2}
  \leq C2^{-k/2}2^{j/2}2^{\theta k/2}\|f_{k,j}(\xi,\tau)\|_{L_{\xi,\tau}^2}
\nonumber \\
  &=& C2^{-(1-\theta)k/2}2^{j/2}\|f_{k,j}(\xi,\tau)\|_{L_{\xi,\tau}^2}
  \leq C2^{-(1-\theta)k/2}2^{j/2}\beta_{k,j}
  \|f_{k,j}(\xi,\tau)\|_{L_{\xi,\tau}^2}.
\nonumber
\end{eqnarray}
  If $j\geq 2k+1$ then we have
$$
  |\tau|\leq |\tau-\omega(\xi)|+|\omega(\xi)|\leq
   C2^j+ C2^{2k}\leq C2^{j},
$$
  so that
\begin{eqnarray}
  &&\|\mathcal{F}^{-1}[|\tau|^{{\theta\over 4}}f_{k,j}(\xi,\tau)]\|_{L_x^\infty L_t^2}
  \leq C2^{-k/2}2^{j/2}2^{\theta j/4}
  \|f_{k,j}(\xi,\tau)\|_{L_{\xi,\tau}^2}
\nonumber \\
  &=& C2^{-(1-\theta)k/2}2^{j/2}2^{\theta(j-2k)/4}
  \|f_{k,j}(\xi,\tau)\|_{L_{\xi,\tau}^2}
  \leq C2^{-(1-\theta)k/2}2^{j/2}\beta_{k,j}
  \|f_{k,j}(\xi,\tau)\|_{L_{\xi,\tau}^2}.
\nonumber
\end{eqnarray}
  Hence
\setcounter{equation}{7}
\begin{eqnarray}
  &&\|D_t^{{\theta\over 4}}\mathcal{F}^{-1}(f_k)\|_{L_x^\infty L_t^2}
  \leq\sum_{j=0}^\infty
  \|\mathcal{F}^{-1}[|\tau|^{{\theta\over 4}}f_{k,j}(\xi,\tau)]\|_{L_x^\infty L_t^2}
\nonumber\\
  &\leq & C2^{-(1-\theta)k/2}\sum_{j=0}^\infty
  2^{j/2}\beta_{k,j}\|f_{k,j}(\xi,\tau)\|_{L_{\xi,\tau}^2}
  =C2^{-(1-\theta)k/2}\|f_k\|_{X_k}.
\end{eqnarray}

  We next assume that $k\geq 100$ and $f_k\in Y_k$. Then $\supp f_k
\subseteq\cup_{j=0}^{k-1}D_{k,j}$, and
\setcounter{equation}{8}
\begin{eqnarray}
  \|f_k\|_{Y_k} &=& 2^{-k/2}\|\mathcal{F}^{-1}[(\tau-\omega(\xi)+i)
  f_k(\xi,\tau)]\|_{L_x^1 L_t^2}
\nonumber \\
  &=& 2^{-k/2}\|F_1^{-1}[(\tau-\omega(\xi)+i)
  f_k(\xi,\tau)]\|_{L_x^1 L_\tau^2}.
\end{eqnarray}
  Let $g_k(x,\tau)=2^{-k/2}F_1^{-1}[(\tau-\omega(\xi)+i)
  f_k(\xi,\tau)]=c2^{-k/2}\displaystyle\int_{-\infty}^\infty
  e^{ix\xi}(\tau-\omega(\xi)+i)f_k(\xi,\tau)d\xi$. Then
  $f_k(\xi,\tau)=2^{k/2}(\tau-\omega(\xi)+i)^{-1}
  \displaystyle\int_{-\infty}^\infty e^{-ix\xi}g_k(x,\tau)dx$
  and, by (3.9), $\|f_k\|_{Y_k}=\|g_k\|_{L_x^1 L_\tau^2}$.
  By the fact that $\supp f_k\subseteq\cup_{j=0}^{k-1}D_{k,j}$ we have
  $f_k(\xi,\tau)=\psi_k(\xi)\eta_0(2^{-k}(\tau-\omega(\xi)))
  f_k(\xi,\tau)$, where $\psi_k(\xi)=\eta_0(2^{-(k+1)}\xi)-
  \eta_0(2^{-(k-2)}\xi)$, so that
$$
  f_k(\xi,\tau)=2^{k/2}\psi_k(\xi)\eta_0(2^{-k}(\tau-\omega(\xi)))
  (\tau-\omega(\xi)+i)^{-1}\int_{-\infty}^\infty
  e^{-iy\xi}g_k(y,\tau)dy.
$$
  Let
$$
  h_k(y,\xi,\tau)=2^{k/2}\psi_k(\xi)\eta_0(2^{-k}(\tau-
\omega(\xi)))(\tau-\omega(\xi)+i)^{-1}e^{-iy\xi}g_k(y,\tau).
\eqno{(3.10)}
$$
  Then the above calculation shows that
$$
  f_k(\xi,\tau)=\int_{-\infty}^\infty h_k(y,\xi,\tau)dy,
\eqno{(3.11)}
$$
  so that
$$
  D_t^{{\theta\over 4}}\mathcal{F}^{-1}(f_k)(x,t)=
  \int_{-\infty}^\infty\!\!\Big(\!\!
  \int_{-\infty}^\infty\!\!\int_{-\infty}^\infty\!\!
  e^{ix\xi}e^{it\tau}|\tau|^{{\theta\over 4}}h_k(y,\xi,\tau)
  d\xi d\tau\Big)dy.
\eqno{(3.12)}
$$
  In what follows we prove that
$$
  \|\int_{-\infty}^\infty\!\!\int_{-\infty}^\infty\!\!
  e^{ix\xi}e^{it\tau}|\tau|^{{\theta\over 4}}h_k(y,\xi,\tau)
  d\xi d\tau\|_{L_x^\infty L_t^2}\leq
  C2^{-(1-\theta)k/2}\|g_k(y,\cdot)\|_2,
\eqno{(3.13)}
$$
  where $C$ is independent of $k$ and $y$. If this inequality is proved,
  then by (3.11) we have
$$
\begin{array}{rl}
  \|D_t^{{\theta\over 4}}\mathcal{F}^{-1}(f_k)\|_{L_x^\infty L_t^2}=
  &\displaystyle\|\int_{-\infty}^\infty\Big(
  \int_{-\infty}^\infty\!\!\int_{-\infty}^\infty\!\!e^{ix\xi}
  e^{it\tau}|\tau|^{{\theta\over 4}}h_k(y,\xi,\tau)
  d\xi d\tau\Big)dy\|_{L_x^\infty L_t^2}
\\ [0.5cm]
  \leq &\displaystyle\int_{-\infty}^\infty\|\int_{-\infty}^\infty
  \!\!\int_{-\infty}^\infty\!\!e^{ix\xi}e^{it\tau}
  |\tau|^{{\theta\over 4}}h_k(y,\xi,\tau)d\xi d\tau\|_{L_x^\infty L_t^2}dy
\\ [0.5cm]
  \leq & C2^{-(1-\theta)k/2}\displaystyle\int_{-\infty}^\infty
  \|g_k(y,\cdot)\|_2 dy
  = C2^{-(1-\theta)k/2}\|g_k\|_{L_x^1 L_\tau^2},
\end{array}
$$
  which, combined with the fact that $\|f_k\|_{Y_k}=\|g_k\|_{L_x^1 L_\tau^2}$,
  yields the following estimate:
$$
  \|D_t^{{\theta\over 4}}\mathcal{F}^{-1}(f_k)\|_{L_x^\infty L_t^2}\leq
  C2^{-(1-\theta)k/2}\|f_k\|_{Y_k}.
\eqno{(3.14)}
$$

  We neglect the parameter $y$ in (3.10) and (3.13). By the Plancherel's
  theorem, (3.13) follows if we prove that
$$
  \Big\|\int_{-\infty}^\infty\!\!e^{ix\xi}|\tau|^{{\theta\over 4}}
  h_k(\xi,\tau)d\xi\Big\|_{L_x^\infty L_\tau^2}\leq
  C2^{-(1-\theta)k/2}\|g_k\|_2.
\eqno{(3.15)}
$$
  To prove this estimate, we first recall that for $k\geq 100$ (see (4.22)
  in \cite{IK}),
$$
  \Big|\int_{-\infty}^\infty\!\!e^{ix\xi}\psi_k(\xi)
  \eta_0(2^{-k}(\tau-\omega(\xi)))(\tau-\omega(\xi)+i)^{-1}d\xi\Big|
  \leq C2^{-k}
\eqno{(3.16)}
$$
  uniformly for $x$ and $\tau$. Next, we note that on the support of $h_k$
  we have $|\xi|\leq C2^k$ and $|\tau-\omega(\xi)|\leq C2^k$, which implies
  that $|\tau|\leq C2^{2k}$. Hence, the left-hand side of (3.15) is dominated
  by
$$
\begin{array}{rl}
  & \displaystyle 2^{k/2}\cdot\sup_{x,\tau}\Big|\int_{-\infty}^\infty\!\!
  e^{ix\xi}\psi_k(\xi)\eta_0(2^{-k}(\tau-\omega(\xi)))
  (\tau-\omega(\xi)+i)^{-1}d\xi\Big|\cdot
  \||\tau|^{{\theta\over 4}}g_k(\tau)\|_{L_{|\tau|\leq C2^{2k}}^2}
\\ [0.3cm]
  \leq & 2^{k/2}\cdot C2^{-k}\cdot C2^{\theta k/2}\|g_k\|_2
  = C2^{-(1-\theta)k/2}\|g_k\|_2,
\end{array}
$$
  as desired.

  By (3.8) and (3.14), we see that (3.5) holds for $k\geq 1$. We now consider
  the case $k=0$. If $f_0\in X_0$ then by (3.2) we have
$$
  D_t^{{\theta\over 4}}\mathcal{F}^{-1}[f_0(\xi,\tau)]=\sum_{j=0}^\infty
  \sum_{l=-\infty}^1 D_t^{{\theta\over 4}}
  \mathcal{F}^{-1}[\phi_0(\xi)\phi_j(\tau)]\ast
  \mathcal{F}^{-1}[f_{0jl}(\xi,\tau)],
$$
  so that
$$
\begin{array}{rl}
  &\|D_t^{{\theta\over 4}}\mathcal{F}^{-1}(f_0)\|_{L_x^\infty L_t^2}\leq
  \displaystyle\sum_{j=0}^\infty\sum_{l=-\infty}^1
  \|D_t^{{\theta\over 4}}\mathcal{F}^{-1}[\phi_0(\xi)\phi_j(\tau)]\|_{L_x^2 L_t^1}
  \|\mathcal{F}^{-1}[f_{0jl}(\xi,\tau)]\|_{L_{x,t}^2}
\\ [0.3cm]
  \leq & C\displaystyle\sum_{j=0}^\infty\sum_{l=-\infty}^1
  \|\phi_0(\xi)\|_{L_\xi^2}\|2^{j(1+{\theta\over 4})}\tilde{\psi_0}(2^jt)]\|_{L_t^1}
  \|f_{0jl}(\xi,\tau)\|_{L_{\xi,\tau}^2}
  \quad (\tilde{\psi_0}=\mathcal{F}^{-1}(\psi_0))
\\ [0.3cm]
  \leq & C\displaystyle\sum_{j=0}^\infty\sum_{l=-\infty}^1
  2^{{\theta j\over 4}}\|f_{0jl}(\xi,\tau)\|_{L_{\xi,\tau}^2}
  \leq C\|f_0\|_{X_0}.
\end{array}
$$
  If $f_0\in Y_0$ then by (3.3) we have
$$
  D_t^{{\theta\over 4}}\mathcal{F}^{-1}[f_0(\xi,\tau)]=\sum_{j=0}^\infty
  D_t^{{\theta\over 4}}\mathcal{F}^{-1}[\phi_0(\xi)\phi_j(\tau)]\ast
  \mathcal{F}^{-1}[f_{0j}(\xi,\tau)],
$$
  so that
$$
\begin{array}{rl}
  &\|D_t^{{\theta\over 4}}\mathcal{F}^{-1}(f_0)\|_{L_x^\infty L_t^2}\leq
  \displaystyle\sum_{j=0}^\infty \|D_t^{{\theta\over 4}}
  \mathcal{F}^{-1}[\phi_0(\xi)\phi_j(\tau)]\|_{L_x^\infty L_t^1}
  \|\mathcal{F}^{-1}[f_{0j}(\xi,\tau)]\|_{L_x^1L_t^2}
\\ [0.3cm]
  \leq & C\displaystyle\sum_{j=0}^\infty
  \|\phi_0(\xi)\|_{L_\xi^1}\|2^{j(1+{\theta\over 4})}
  \tilde{\psi_0}(2^jt)]\|_{L_t^1}
  \|\mathcal{F}^{-1}[\eta_j(\tau)f_0(\xi,\tau)]\|_{L_x^1L_t^2}
 \qquad (\tilde{\psi_0}=\mathcal{F}^{-1}(\psi_0))
\\ [0.3cm]
  \leq & C\displaystyle\sum_{j=0}^\infty 2^{{\theta j\over 4}}
  \|\mathcal{F}^{-1}[\eta_j(\tau)f_0(\xi,\tau)]\|_{L_x^1L_t^2}
  \leq C\|f_0\|_{Y_0}.
\end{array}
$$
  Hence (3.5) also holds for $k=0$. The proof is complete. $\quad\Box$
\medskip

  {\bf Lemma 3.4}\ \ {\it If $w\in F^0$ then for any $0\leq\theta<1$ we have
  $D_t^{{\theta\over 4}}w\in L_x^\infty L_t^2$, and}
$$
  \|D_t^{{\theta\over 4}}w\|_{L_x^\infty L_t^2}\leq C\|w\|_{F^0}.
\eqno{(3.17)}
$$

  {\it Proof}:\ \ Let $f_k(\xi,\tau)=\eta_k(\xi)\widetilde{w}(\xi,\tau)$,
  $k=0,1,2,\cdots$, where $\widetilde{w}=F(w)$. Then $w\in F^0$ implies
  that $(I-\partial_\tau^2)f_k\in Z_k$, $k=0,1,2,\cdots$, and
$$
 \|w\|_{F^0}=\Big(\sum_{k=0}^\infty\|(I-\partial_\tau^2)f_k\|_{Z_k}^2
  \Big)^{{1\over 2}}<\infty.
$$
  Since $\widetilde{w}(\xi,\tau)=\displaystyle\sum_{k=0}^\infty
  f_k(\xi,\tau)$, for any $0\leq\theta<1$ we have
$$
  D_t^{{\theta\over 4}}[(1+t^2)w(x,t)]=\sum_{k=0}^\infty
  D_t^{{\theta\over 4}}\mathcal{F}^{-1}[(I-\partial_\tau^2)f_k(\xi,\tau)].
$$
  Hence, by Lemma 3.3 we have
\setcounter{equation}{17}
\begin{eqnarray}
  &&\|D_t^{{\theta\over 4}}[(1+t^2)w(x,t)]\|_{L_x^\infty L_t^2}
  \leq\displaystyle\sum_{k=0}^\infty\|D_t^{{\theta\over 4}}
  \mathcal{F}^{-1}[(I-\partial_\tau^2)f_k(\xi,\tau)]\|_{L_x^\infty L_t^2}
\nonumber \\
  &\leq & C\displaystyle\sum_{k=0}^\infty 2^{-(1-\theta)k/2}
  \|(I-\partial_\tau^2)f_k(\xi,\tau)\|_{Z_k}
\nonumber \\
  &\leq & C\displaystyle
  \Big(\sum_{k=0}^\infty\|(I-\partial_\tau^2)
  f_k(\xi,\tau)\|_{Z_k}^2\Big)^{{1\over 2}}
  =C\|w\|_{F^0}.\quad
\end{eqnarray}
  Since $w(x,t)=(1+t^2)^{-1}\cdot (1+t^2)w(x,t)$, by Theorem A.12 in
  \cite{KPV2} we have
$$
\begin{array}{rl}
  &\|D_t^{{\theta\over 4}}w-(1+t^2)^{-1}D_t^{{\theta\over 4}}
  [(1+t^2)w(x,t)]-D_t^{{\theta\over 4}}(1+t^2)^{-1}\cdot
  (1+t^2)w(x,t)\|_{L_x^\infty L_t^2}\\
  \leq &C\|(1+t^2)^{-1}\|_\infty\|D_t^{{\theta\over 4}}
  [(1+t^2)w(x,t)]\|_{L_x^\infty L_t^2}
\end{array}
$$
  From this estimate and (3.18), we see that (3.17) follows. $\quad\Box$
\medskip

  {\bf Lemma 3.5}\ \ {\it Let $f_k\in Z_k$, $k\geq 0$. Then for any
  admissible pair $(p,q)$ we have}
$$
  \|\mathcal{F}^{-1}(f_k)\|_{L_t^q L_x^p}\leq C(p,q)\|f_k\|_{Z_k}.
\eqno{(3.19)}
$$

  {\it Proof}:\ \ Assume first that $k\geq 1$ and $f_k\in X_k$. Let
  $f_{k,j}(\xi,\tau)=\eta_j(\tau-\omega(\xi))f_k(\xi,\tau)$, $j\in
  \mathbb{Z}$, $j\geq 0$. Then $\supp f_{k,j}\subseteq D_{k,j}$, and
  (3.6), (3.7) hold. Let $f_{k,j}^\#(\xi,\tau)=f_{k,j}(\xi,\tau+\omega(\xi))$.
  Then $\supp f_{k,j}^\#\subseteq I_k\times\tilde{I}_j$. We have
$$
\begin{array}{rl}
  \mathcal{F}^{-1}(f_{k,j})=&c^2\displaystyle\int_{-\infty}^\infty
  \int_{-\infty}^\infty e^{it\tau}e^{ix\xi}f_{k,j}(\xi,\tau)d\xi d\tau
\\ [0.3cm]
  =&c^2\displaystyle\int_{-\infty}^\infty\int_{-\infty}^\infty
   e^{it\tau}e^{ix\xi}e^{it\omega(\xi)}
   f_{k,j}(\xi,\tau+\omega(\xi))d\xi d\tau
\\ [0.3cm]
  =&c^2\displaystyle\int_{\tilde{I}_j} e^{it\tau}\Big(
  \int_{-\infty}^\infty e^{ix\xi}e^{it\omega(\xi)}
  f_{k,j}^\#(\xi,\tau)d\xi\Big)d\tau
\end{array}
$$
  Let $g_{k,j}^\tau(x)=c\displaystyle\int_{-\infty}^\infty e^{ix\xi}
  f_{k,j}^\#(\xi,\tau)d\xi=F_1^{-1}(f_{k,j}^\#(\cdot,\tau))$. Then
  $f_{k,j}^\#(\xi,\tau)=F_1(g_{k,j}^\tau)$, so that
$$
  \mathcal{F}^{-1}(f_{k,j})=c\int_{\tilde{I}_j} e^{it\tau}F_1^{-1}
  \Big(e^{it\omega(\xi)}F_1(g_{k,j}^\tau)\Big)d\tau
  =c\int_{\tilde{I}_j} e^{it\tau}W(t)g_{k,j}^\tau(x)d\tau.
$$
  It follows that
\setcounter{equation}{19}
\begin{eqnarray}
  \|\mathcal{F}^{-1}(f_{k,j})\|_{L_t^q L_x^p}
  &\leq & c\displaystyle\int_{\tilde{I}_j}
  \|W(t)g_{k,j}^\tau(x)\|_{L_t^q L_x^p}d\tau
  \leq C\displaystyle\int_{\tilde{I}_j}
  \|g_{k,j}^\tau(x)\|_2 d\tau
\nonumber\\
  && \qquad\qquad\qquad\qquad\qquad
  (\mbox{by Strichartz estimate})
\nonumber\\
  &=& C\displaystyle\int_{\tilde{I}_j}
  \Big[\int_{-\infty}^\infty |f_{k,j}^\#(\xi,\tau)|^2
  d\xi\Big]^{{1\over 2}}d\tau \leq
  C2^{j/2}\|f_{k,j}\|_{L^2_{\xi,\tau}}.
\end{eqnarray}
  By (3.6), (3.7) and (3.20) we have
$$
  \|\mathcal{F}^{-1}(f_{k})\|_{L_t^q L_x^p}\leq
  \sum_{j=0}^\infty\|\mathcal{F}^{-1}(f_{k,j})\|_{L_t^q L_x^p}
  \leq C\sum_{j=0}^\infty2^{j/2}\|f_{k,j}\|_{L^2_{\xi,\tau}}
  \leq C\|f_{k}\|_{X_k}.
\eqno{(3.21)}
$$

  Next assume that $k\geq 1$ and $f_k\in Y_k$. Then $\supp f_k\subseteq
  \cup_{j=0}^{k-1}D_{k,j}$, and (3.9) holds. Let $g_k(x,\tau)$ and
  $h_k(y,\xi,\tau)$ be as in the proof of Lemma 3.3. In what follows
  we prove that
$$
  \|\int_{-\infty}^\infty\!\!\int_{-\infty}^\infty\!\!
  e^{ix\xi}e^{it\tau}h_k(y,\xi,\tau)d\xi d\tau\|_{L_t^q L_x^p}
  \leq C\|g_k(y,\cdot)\|_2,
\eqno{(3.22)}
$$
  where $C$ is independent of $k$ and $y$. If this inequality is proved,
  then by (3.11) we have
$$
\begin{array}{rl}
  \|\mathcal{F}^{-1}(f_{k})\|_{L_t^q L_x^p}=
  &\displaystyle\|\int_{-\infty}^\infty\Big(
  \int_{-\infty}^\infty\!\!\int_{-\infty}^\infty\!\!e^{i(x-y)\xi}
  e^{it\tau}h_k(y,\xi,\tau)d\xi d\tau\Big)dy\|_{L_t^q L_x^p}
\\ [0.5cm]
  \leq &\displaystyle\int_{-\infty}^\infty\|\int_{-\infty}^\infty
  \!\!\int_{-\infty}^\infty\!\!e^{i(x-y)\xi}e^{it\tau}
  h_k(y,\xi,\tau)d\xi d\tau\|_{L_t^q L_x^p}dy
\\ [0.5cm]
  \leq & C\displaystyle\int_{-\infty}^\infty
  \|g_k(y,\cdot)\|_2 dy
  =C\|g_k\|_{L_x^1 L_\tau^2},
\end{array}
$$
  which, combined with the fact that $\|g_k\|_{L_x^1 L_\tau^2}=\|f_k\|_{Y_k}$,
  gives the desired assertion.

  We neglect the parameter $y$ in (3.10) and (3.22). Since $k\geq 100$ and
  $|\xi|\in [2^{k-2},2^{k+2}]$, we may assume that the function $g_k=
  g_k(\tau)$ in (3.10) is supported in the set $\{\tau:|\tau|\in [2^{2k-10},
  2^{2k+10}]\}$. Let $g_k^+=g_k\cdot\chi_{[0,\infty)}$, $g_k^-=g_k\cdot
  \chi_{(-\infty,0]}$, and define the corresponding function $h_k^+$ and
  $h_k^-$ as in (3.10). By symmetry, it suffices to prove (2.22) for $h_k^+$,
  which is supported in $\{(\xi,\tau):\xi\in [-2^{k-2},-2^{k+2}],\tau\in
  [2^{2k-10},2^{2k+10}]\}$. Since $\omega(\xi)=-\xi|\xi|$, we have
  $\tau-\omega(\xi)=\tau-\xi^2$ on the support of $h_k^+$, and $h_k^+(\xi,
  \tau)=0$ unless $|\sqrt{\tau}+\xi|\leq C$. Let
$$
  {h'}_k^+(\xi,\tau)=2^{k/2}\psi_k(-\sqrt{\tau})\eta_0(\sqrt{\tau}+\xi)
  [\tau-\xi^2+(\sqrt{\tau}+\xi)+i\sqrt{\tau}2^{-k}]^{-1}g_k^+(\tau).
$$
  By the argument in Lines 25--28 in Page 762 of \cite{IK}, we know that
$$
  \|h_k^+-{h'}_k^+\|_{X_k}\leq C\|g_k^+\|_2.
$$
  Hence, by (3.21) we have
$$
  \|\mathcal{F}^{-1}(h_k^+-{h'}_k^+)\|_{L_t^q L_x^p}\leq C\|g_k^+\|_2.
\eqno{(3.23)}
$$
  To estimate $\|\mathcal{F}^{-1}({h'}_k^+)\|_{L_t^q L_x^p}$, we make the
  change of variables $\xi\to\xi'$ by letting $\xi=\xi'-\sqrt{\tau}$. Then
  we have
$$
\begin{array}{rl}
  \mathcal{F}^{-1}({h'}_k^+)(x,t)=&
 \displaystyle 2^{k/2}\int_{-\infty}^\infty e^{it\tau}e^{-ix\sqrt{\tau}}
  \psi_k(-\sqrt{\tau})g_k^+(\tau)(2\sqrt{\tau})^{-1}d\tau
\\ [0.3cm]
 &\quad\times\displaystyle \int_{-\infty}^\infty e^{ix\xi'}
  \eta_0(\xi')(\xi'+i2^{-k-1})^{-1}d\xi'.
\end{array}
$$
  It can be easily seen that the second integral is bounded by a constant
  independent of $x$ and $k$. Next we compute
$$
\begin{array}{rll}
  &&\displaystyle \|\int_{-\infty}^\infty e^{it\tau}e^{-ix\sqrt{\tau}}
  \psi_k(-\sqrt{\tau})g_k^+(\tau)(2\sqrt{\tau})^{-1}d\tau\|_{L_t^q L_x^p}
\\ [0.3cm]
  &=&\displaystyle \|\int_{-\infty}^\infty e^{it\xi^2}e^{-ix\xi}
  \psi_k(-\xi)g_k^+(\xi^2)d\xi\|_{L_t^q L_x^p}
  \quad (\mbox{by letting}\;\;\tau=\xi^2)
\\ [0.3cm]
  &\leq & C\displaystyle\Big[\int_{-\infty}^\infty
  |\psi_k(-\xi)g_k^+(\xi^2)|^2d\xi\Big]^{1/2}
  \quad (\mbox{by using Strichartz and Plancherel})
\\ [0.3cm]
  &=& C\displaystyle\Big[\int_{-\infty}^\infty
  |\psi_k(-\sqrt{\tau})g_k^+(\tau)|^2(2\sqrt{\tau})^{-1}d\tau\Big]^{1/2}
\\ [0.3cm]
  &\leq & C2^{-k/2}\|g_k^+\|_2
  \qquad\qquad\qquad\qquad (\mbox{because}\;\;\sqrt{\tau}\sim 2^k)
\end{array}
$$
  Hence
$$
  \|\mathcal{F}^{-1}({h'}_k^+)\|_{L_t^q L_x^p}\leq C\|g_k^+\|_2.
$$
  Combining this estimate with (3.23), we see that the desired assertion
  follows.

  From the above deduction we see that (3.19) holds for $k\geq 1$. For the
  case $k=0$, the argument is similar to that in the proof of  Lemma 3.1.
  Indeed, if $f_0\in X_0$ then from (3.2) we have
$$
\begin{array}{rl}
  \|\mathcal{F}^{-1}(f_0)\|_{L_t^q L_x^p}\leq
  &\displaystyle\sum_{j=0}^\infty\sum_{l=-\infty}^1
  \|\mathcal{F}^{-1}[\phi_0(\xi)\phi_j(\tau)]\|_{L_t^s L_x^r}
  \|\mathcal{F}^{-1}[f_{0jl}(\xi,\tau)]\|_{L_{x,t}^2}
\\ [0.3cm]
&\quad\quad\quad\quad\quad\quad\quad\quad\quad\quad
   (1/r=1/p+1/2, \quad 1/s=1/q+1/2)
\\ [0.3cm]
  \leq & C\displaystyle\sum_{j=0}^\infty\sum_{l=-\infty}^1
  \|\mathcal{F}^{-1}(\phi_0)\|_{L_x^r}\|2^j\tilde{\psi_0}(2^jt)]\|_{L_t^s}
  \|f_{0jl}(\xi,\tau)\|_{L_{\xi,\tau}^2}
 \qquad (\tilde{\psi_0}=\mathcal{F}^{-1}(\psi_0))
\\ [0.3cm]
  \leq & C\displaystyle\sum_{j=0}^\infty\sum_{l=-\infty}^1
  2^{j(1-1/s)}\|\eta_j(\tau)\chi_l(\xi)f_0(\xi,\tau)\|_{L_{\xi,\tau}^2}
  \leq C\|f_0\|_{X_0}.
\end{array}
$$
  If $f_0\in Y_0$ then from (3.3) we have
$$
\begin{array}{rl}
  \|\mathcal{F}^{-1}(f_0)\|_{L_t^q L_x^p}\leq
  &\displaystyle\sum_{j=0}^\infty
  \|\mathcal{F}^{-1}[\phi_0(\xi)\phi_j(\tau)]\|_{L_t^s L_x^p}
  \|\mathcal{F}^{-1}[f_{0j}(\xi,\tau)]\|_{L_t^2L_x^1}\quad
   (1/s=1/q+1/2)
\\ [0.3cm]
  \leq & C\displaystyle\sum_{j=0}^\infty
  \|\mathcal{F}^{-1}(\phi_0)\|_{L_x^p}\|2^j\tilde{\psi_0}(2^jt)]\|_{L_t^s}
  \|\mathcal{F}^{-1}[\eta_j(\tau)f_0(\xi,\tau)]\|_{L_x^1L_t^2}
 \quad (\tilde{\psi_0}=\mathcal{F}^{-1}(\psi_0))
\\ [0.3cm]
  \leq & C\displaystyle\sum_{j=0}^\infty
  2^{j(1-1/s)}\|\mathcal{F}^{-1}[\eta_j(\tau)f_0(\xi,\tau)]\|_{L_x^1L_t^2}
  \leq C\|f_0\|_{Y_0}.
\end{array}
$$
  Hence the desired assertion also holds for $k=0$. $\quad\Box$
\medskip

  Using the above lemma and a similar argument to the one in the proof of Lemma
  3.2 we have
\medskip

  {\bf Lemma 3.6}\ \ {\it Assume that $w\in F^0$. Then for any admissible pair
  $(p,q)$ we have $w\in L_t^q(\mathbb{R}, L_x^p(\mathbb{R}))$, and}
$$
  \|w\|_{L_t^q L_x^p}\leq C_{pq}\|w\|_{F^0}.
\eqno{(3.24)}
$$

  Since (6,6) is an admissible pair, by the above lemma we have
\medskip

  {\bf Corollary 3.7}\ \ {\it Assume that $w\in F^0$. Then $w\in
  L^6(\mathbb{R}^2)$, and}
$$
  \|w\|_{L_{x,t}^6}\leq C\|w\|_{F^0}.
\eqno{(3.25)}
$$

  Using the expression (2.5) and Lemma 3.6 we have
\medskip

  {\bf Corollary 3.8}\ \ {\it Let $\phi\in L^2(\mathbb{R})$ and let $u$ be the
  global solution of the problem $(1.1)$ ensured by Theorem 2.2. Then for any
  $T>0$ and any admissible pair $(p,q)$ we have $u\in L^q([-T,T],
  L_x^p(\mathbb{R}))$. Moreover, the mapping $\phi\to u$ from $L^2(\mathbb{R})$
  to $L^q([-T,T],L_x^p(\mathbb{R}))$ defined in this way is continuous, and
  there exists corresponding function $c_{pqT}:[0,\infty)\to [0,\infty)$ such
  that}
$$
  \|u\|_{L_T^q L_x^p}\leq c_{pqT}(\|\phi\|_2).
\eqno{(3.26)}
$$

  {\it Proof}:\ \ Choose $M\geq\|\phi\|_2$ and fix it. With $M$
  fixed in this way, let $\delta$ be as in Corollary 3.4. By dividing the
  interval $[-T,T]$ into subintervals $[-\delta,\delta]$, $\pm [\delta,2\delta]$,
  $\pm [2\delta,3\delta]$, $\cdots$, $\pm [(N\!-\!1)\delta,T]$, where $N$ is
  the smallest integer such that $T\leq N\delta$, and using the $L^2$-conservation
  law, we only need to prove the assertion holds for $T=\delta$. For $T=\delta$
  the expression (2.5) holds, from which the desired assertion easily follows.
  Indeed, from (2.2) it is clear that for any admissible pair $(p,q)$ we have
  $u_0\in L^q([-T,T],L_x^p(\mathbb{R}))$,
$$
  \|u_0\|_{L_T^q L_x^p}\leq c_{pqT}(\|\phi\|_2),
\eqno{(3.27)}
$$
  and the mapping $\phi\to u_0$ from $L^2(\mathbb{R})$ to $L^q([-T,T],
  L_x^p(\mathbb{R}))$ is continuous. Secondly, since $U_0$ is real, we see that
  $e^{\pm iU_0}$ are uniformly bounded, and it is clear that the mappings
  $\phi\to e^{\pm iU_0}$ from $L^2(\mathbb{R})$ to $L^\infty(\mathbb{R}\times
  [-T,T])$ are continuous. Finally, by Theorem 2.1, Lemma 3.6, and the continuity
  of the mapping $\phi\to (e^{iU_0(\cdot,0)}\phi_{+{\rm high}},e^{-iU_0(\cdot,0)}
  \phi_{-{\rm high}},0)$ from $L^2(\mathbb{R})$ to $(\widetilde{H}^0)^3$, we
  see that for any admissible pair $(q,p)$ we have $w_+,w_-,w_0\in L^q([-T,T],
  L_x^p(\mathbb{R}))$,
$$
  \|w_+\|_{L_T^q L_x^p}+\|w_-\|_{L_T^q L_x^p}+\|w_0\|_{L_T^q L_x^p}\leq
  c_{pqT}(\|\phi\|_2),
\eqno{(3.28)}
$$
  and the mapping $\phi\to (w_+,w_-,w_0)$ from $L^2(\mathbb{R})$ to $[L^q([-T,T],
  L_x^p(\mathbb{R}))]^3$ is continuous. By (2.5), (3.27), (3.28) and the uniform
  boundedness of $e^{\pm iU_0}$ we have
$$
  \|u\|_{L_T^q L_x^p}\leq\|w_+\|_{L_T^q L_x^p}+\|w_-\|_{L_T^q L_x^p}
  +\|w_0\|_{L_T^q L_x^p}+\|u_0\|_{L_T^q L_x^p}\leq
  c_{pqT}(\|\phi\|_2).
$$
  Hence $u\in L^q([-T,T],L_x^p(\mathbb{R}))$ and (2.26) holds. Moreover, the
  above argument also shows that the mapping $\phi\to u$ from $L^2(\mathbb{R})$
  to $L^q([-T,T],L_x^p(\mathbb{R}))$ is continuous. The proof is complete.
$\quad\Box$

\section{The proof of Theorem 1.1}

  Since we are not clear if the space $F^0$ in which uniqueness of the solution
  of (2.4) is ensured is reflexive, we cannot use functional analysis to get
  the assertion that any bounded sequence in $F^0$ has a weakly convergent
  subsequence. To overcome this difficulty, we shall appeal to the following
  preliminary result:
\medskip

  {\bf Lemma 4.1}\ \ {\em Let $w_n\in F^0\cap L^2(\mathbb{R}^2)$, $n=1,2,\cdots$.
  Assume that $\|w_n\|_{F^0}\leq M$ for all $n\in\mathbb{N}$ and some $M>0$, and
  there exists $T>0$ such that $w_n(t)=0$ for all $|t|\geq T$ and $n\in\mathbb{N}$.
  Assume further that as $n\to\infty$, $w_n\to w$ weakly in $L^2(\mathbb{R}^2)$.
  Then $w\in F^0$, and $\|w\|_{F^0}\leq M$, or more precisely,}
$$
  \|w\|_{F^0}\leq\liminf_{n\to\infty}\|w_n\|_{F^0}.
\eqno{(4.1)}
$$

  {\it Proof}:\ \ We fulfill the proof in three steps.
\medskip

  {\it Step 1}:\ \ We first prove that similar results hold for the spaces $Y_k$
  and $X_k$. That is, taking $Y_k$ as an example and assuming that $f_n\in Y_k
  \cap L^2(\mathbb{R}^2)$, $n=1,2,\cdots$, $\|f_n\|_{Y_k}\leq M$ for some $M>0$
  and all $n\in\mathbb{N}$, and as $n\to\infty$, $f_n\to f$ weakly in
  $L^2(\mathbb{R}^2)$, we have that $f\in Y_k$, and $\|f\|_{Y_k}\leq M$. Note
  that if this assertion is proved, then it follows immediately that
$$
  \|f\|_{Y_k}\leq\liminf_{n\to\infty}\|f_n\|_{Y_k}.
$$

  Consider first the case $k\geq 1$. Let $\psi\in C^\infty_0(\mathbb{R}^2)$ be
  such that
$$
  0\leq\psi\leq 1, \quad \mbox{and} \quad
  \psi(x,t)=1\;\; \mbox{for}\;\; |x|\leq 1,\;\; |t|\leq 1.
$$
  Let $\psi_R(x,t)=\psi(x/R,t/R)$, $R>1$. Since $\|\psi_R\|_{L_{x,t}^\infty}=1$,
  we have, for any $R>1$ and $n\in\mathbb{N}$,
$$
  \|\psi_R \mathcal{F}^{-1}[(\tau-\omega(\xi)+i)f_n(\xi,\tau)]\|_{L_x^1 L_t^2}
  \leq \|\psi_R\|_{L_{x,t}^\infty}\cdot 2^{k/2}\|f_n\|_{Y_k}\leq 2^{k/2}M.
\eqno{(4.2)}
$$
  We first assume that as $n\to\infty$, $f_n\to f$ strongly in $L^2(\mathbb{R}^2)$.
  Let $\varphi_k\in C_0^\infty(\mathbb{R}^2)$ be such that $\varphi_k(\xi,\tau)=1$
  for $(\xi,\tau)\in \cup_{j=1}^{k-1}D_{k,j}$. Since $\supp f_n\subseteq
  \cup_{j=1}^{k-1}D_{k,j}$ for all $n\in\mathbb{N}$, for any $m,n\in\mathbb{N}$
  we have
$$
\begin{array}{rcl}
  &&\|\psi_R \mathcal{F}^{-1}[(\tau-\omega(\xi)+i)(f_n-f_m)]\|_{L_x^1 L_t^2}
\\ [0.3cm]
  &\leq &\|\psi_R\|_{L_x^2 L_t^\infty}
  \|\mathcal{F}^{-1}[\varphi_k(\xi,\tau)(\tau-\omega(\xi)+i)(f_n-f_m)]\|_{L_{x,t}^2}
\\ [0.3cm]
  &\leq &\|\psi_R\|_{L_x^2 L_t^\infty}
  \|\mathcal{F}^{-1}[\varphi_k(\xi,\tau)(\tau-\omega(\xi)+i)]\|_{L_{x,t}^1}
  \|f_n-f_m\|_{L_{x,t}^2}.
\end{array}
$$
  From this we see that for any $R>1$, $\psi_R \mathcal{F}^{-1}[(\tau-\omega(\xi)+i)
  f_n(\xi,\tau)]$ is convergent in $L_x^1 L_t^2$, so that $\psi_R \mathcal{F}^{-1}
  [(\tau-\omega(\xi)+i)f(\xi,\tau)]\in L_x^1 L_t^2$ and, by letting $n\to\infty$ in
  (4.2) we get, for any $R>0$,
$$
  \|\psi_R \mathcal{F}^{-1}[(\tau-\omega(\xi)+i)f(\xi,\tau)]\|_{L_x^1 L_t^2}
  \leq 2^{k/2}M.
\eqno{(4.3)}
$$
  Clearly, $\mathcal{F}^{-1}[(\tau-\omega(\xi)+i)f(\xi,\tau)$ is a measurable
  function and, as $R\to\infty$, $\psi_R \mathcal{F}^{-1}[(\tau-\omega(\xi)+i)
  f(\xi,\tau)]$ pointwisely converges to $\mathcal{F}^{-1}[(\tau-\omega(\xi)+i)
  f(\xi,\tau)$. Hence, by letting $R\to\infty$ in $(2.3)$ and using Fatou's
  lemma we get
$$
  \|\mathcal{F}^{-1}[(\tau-\omega(\xi)+i)f(\xi,\tau)]\|_{L_x^1 L_t^2}
  \leq 2^{k/2}M,
$$
  so that $f\in Y_k$ and $\|f\|_{Y_k}\leq M$. We next assume that as $n\to\infty$,
  $f_n\to f$ weakly in $L^2(\mathbb{R}^2)$. By a well-known theorem in functional
  analysis, we know that there is another sequence $f_n'$, $n=1,2,\cdots$, with
  each $f_n'$ being a convex combination of finite elements in $\{f_n\}$, such
  that as $n\to\infty$, $f_n'\to f$ strongly in $L^2(\mathbb{R}^2)$. Clearly,
  $\|f_n'\|_{Y_k}\leq M$ for all $n\in\mathbb{N}$. Hence, by the assertion we
  have just proved it follows that $f\in Y_k$ and $\|f\|_{Y_k}\leq M$. This
  proves the desired assertion for the case $k\geq 1$.

  Consider next the case $k=0$. For any $N\in\mathbb{N}$ we have
$$
  \sum_{j=0}^N 2^j\|\mathcal{F}^{-1}[\eta_j(\tau)f_n(\xi,\tau)\|_{L_x^1 L_t^2}
  \leq\|f_n\|_{Y_0}\leq M, \quad n=1,2,\cdots.
$$
  Since for every $0\leq j\leq N$, $\eta_j(\tau)f_n(\xi,\tau)$ have supports
  contained in a common compact set, the argument for the case $k\geq 1$ applies
  to the sequences $\{\eta_j(\tau)f_n(\xi,\tau)\}$, $j=0,1,\cdots, N$, so that
$$
  \|\mathcal{F}^{-1}[\eta_j(\tau)f(\xi,\tau)\|_{L_x^1 L_t^2}\leq
  \liminf_{n\to\infty}\|\mathcal{F}^{-1}[\eta_j(\tau)f_n(\xi,\tau)\|_{L_x^1 L_t^2},
  \quad j=0,1,\cdots, N.
$$
  Hence
$$
  \sum_{j=0}^N 2^j\|\mathcal{F}^{-1}[\eta_j(\tau)f(\xi,\tau)\|_{L_x^1 L_t^2}\leq
  \liminf_{n\to\infty}\sum_{j=0}^N 2^j
  \|\mathcal{F}^{-1}[\eta_j(\tau)f_n(\xi,\tau)\|_{L_x^1 L_t^2}\leq M.
$$
  By the arbitrariness of $N$, we conclude that $f\in Y_0$ and $\|f\|_{Y_0}\leq
  M$, as desired.

  The proof for $X_k$ ($k\geq 0$) follows from a similar argument as in the
  proof for $Y_0$.
\medskip

  {\it Step 2}:\ \ We next prove that a similar result holds for the space $Z_k$,
  namely, assuming that $f_n\in Z_k\cap L^2(\mathbb{R}^2)$, $n=1,2,\cdots$,
  $\|f_n\|_{Z_k}\leq M$ for some $M>0$ and all $n\in\mathbb{N}$, and as $n\to
  \infty$, $f_n\to f$ weakly in $L^2(\mathbb{R}^2)$, then we have $f\in Z_k$, and
  $\|f\|_{Z_k}\leq M$, or more precisely,
$$
  \|f\|_{Z_k}\leq\liminf_{n\to\infty}\|f_n\|_{Z_k}.
$$

  To prove this assertion, we only need to prove that $f\in Z_k$ and, for any
  $\epsln>0$, we have $\|f\|_{Z_k}\leq M+\epsln$. Assume that either $k\geq 100$
  or $k=0$ (the case $1\leq k\leq 99$ is obvious). Given $\epsln>0$, by the
  definition of $Z_k$ we can find for each $n\in\mathbb{N}$ two functions $g_n$
  and $h_n$, $g_n\in X_k$, $h_n\in Y_k$, and
$$
  \|g_n\|_{X_k}+\|h_n\|_{Y_k}\leq\|f_n\|_{Z_k}+\epsln\leq M+\epsln.
\eqno{(4.4)}
$$
  Let $\varphi_k$ be as before. Then we have
$$
  h_n=\varphi_k h_n=\varphi_k f_n-\varphi_k g_n.
\eqno{(4.5)}
$$
  From the definition of $X_k$ and the fact that $g_n\in X_k$ it can be easily
  seen that $\varphi_k g_n\in L^2(\mathbb{R}^2)$, and there exists constant
  $C_k>0$ such that
$$
  \|\varphi_k g_n\|_{L^2(\mathbb{R}^2)}\leq C_k\|g_n\|_{X_k}
  \leq C_k (M+\epsln).
$$
  Hence, there exists a subsequence of $\{g_n\}$, for simplicity of the
  notation we assume that this subsequence is the whole sequence $\{g_n\}$,
  and a function $h_0\in L^2(\mathbb{R}^2)$, such that as $n\to\infty$,
  $\varphi_k g_n\to h_0$ weakly in $L^2(\mathbb{R}^2)$. Let $h=\varphi_k f
  -h_0$. Then $h\in L^2(\mathbb{R}^2)$, and by $(2.5)$ and the fact that
  $f_n\to f$ weakly in $L^2(\mathbb{R}^2)$ we see that $h_n\to h$ weakly in
  $L^2(\mathbb{R}^2)$. Since $\|h_n\|_{Y_k}\leq M+\epsln$, $n=1,2,\cdots$,
  by using the assertion in Step 1 we conclude that $h\in Y_k$, and
$$
  \|h\|_{Y_k}\leq\liminf_{n\to\infty}\|h_n\|_{Y_k}<\infty.
\eqno{(4.6)}
$$
  Now, since both $f_n\to f$ and $h_n\to h$ weakly in $L^2(\mathbb{R}^2)$,
  it follows that $g_n\to g\equiv f-h$ weakly in $L^2(\mathbb{R}^2)$, which
  further implies that for any $j\in\mathbb{Z}\cap [0,\infty)$, $\eta_j(\tau-
  \omega(\xi))g_n(\xi,\tau)\to\eta_j(\tau-\omega(\xi))g(\xi,\tau)$ weakly in
  $L^2(\mathbb{R}^2)$. Since for any $N\in\mathbb{N}$ we have
$$
  \sum_{j=0}^N 2^{j/2}\beta_{k,j}\|\eta_j(\tau-\omega(\xi))
  g_n(\xi,\tau)\|_{L_{\xi,\tau}^2}\leq\|g_n\|_{X_k}\leq M+\epsln,
  \quad n=1,2,\cdots,
$$
  letting $n\to\infty$ we get
$$
\begin{array}{rl}
  \displaystyle\sum_{j=0}^N 2^{j/2}\beta_{k,j}\|\eta_j(\tau-\omega(\xi))
    g(\xi,\tau)\|_{L_{\xi,\tau}^2}
  \leq &\displaystyle\liminf_{n\to\infty}\sum_{j=0}^N 2^{j/2}\beta_{k,j}
  \|\eta_j(\tau-\omega(\xi))g_n(\xi,\tau)\|_{L_{\xi,\tau}^2}
\\ [0.5cm]
  \leq &\displaystyle\liminf_{n\to\infty}\|g_n\|_{X_k} <\infty.
\end{array}
$$
  By arbitrariness of $N$ we conclude that $g\in X_k$, and
$$
  \|g\|_{X_k}\leq\liminf_{n\to\infty}\|g_n\|_{X_k}<\infty.
\eqno{(4.7)}
$$
  Hence, $f=g+h\in Z_k$, and by (4.4), (4.6) and (4.7) we have
$$
  \|f\|_{Z_k}\leq\|g\|_{X_k}+\|h\|_{Y_k}\leq\liminf_{n\to\infty}
  (\|g_n\|_{X_k}+\|h_n\|_{Y_k})\leq M+\epsln.
$$
  This proves the desired assertion.
\medskip

  {\it Step 3}: We now arrive at the last step of the proof of Lemma 4.1. Let
  $f_{nk}(\xi,\tau)=\eta_k(\xi)(I-\partial_\tau^2)\tilde{w}_n(\xi,\tau)$,
  $k=0,1,2,\cdots$, and $f_k(\xi,\tau)=\eta_k(\xi)(I-\partial_\tau^2)
  \tilde{w}(\xi,\tau)$, where $\tilde{w}_n=F(w_n)$ and $\tilde{w}=F(w)$. Then
$$
  \|w_n\|_{F^0}=\Big(\sum_{k=0}^\infty
  \|f_{nk}\|_{Z_k}^2\Big)^{\frac{1}{2}}\leq M,
  \quad n=1,2,\cdots,
$$
  so that for any $N\in\mathbb{N}$ we have
$$
  \sum_{k=0}^N\|f_{nk}\|_{Z_k}^2\leq M^2,
  \quad n=1,2,\cdots,
\eqno{(4.8)}
$$
  Since $w_n$ weakly converges to $w$ in $L^2(\mathbb{R}^2)$ and $w_n(t)=0$ for
  $|t|\geq T$, we have that also $(1+t^2)w_n$ weakly converges to $(1+t^2)w$ in
  $L^2(\mathbb{R}^2)$. By the Parseval formula
$$
  \int\!\!\int_{\mathbb{R}^2}\!\!\tilde{f}(\xi,\tau)
  \varphi(\xi,\tau)d\xi d\tau=
  \int\!\!\int_{\mathbb{R}^2}\!\!f(x,t)\tilde{\varphi}(x,t)dx dt, \quad
  f,\varphi\in L^2(\mathbb{R}^2),
$$
  it follows immediately that $(I-\partial_\tau^2)\tilde{w}_n(\xi,\tau)$
  weakly converges to $(I-\partial_\tau^2)\tilde{w}(\xi,\tau)$ in
  $L^2(\mathbb{R}^2)$, which further implies that for every $k\in\mathbb{Z}\cap
  [0,\infty)$, $f_{nk}$ weakly converges to $f_k$ in $L^2(\mathbb{R}^2)$. Hence,
  by the assertion we proved in Step 2 and $(2.7)$ we get
$$
  \sum_{k=0}^N\|f_k\|_{Z_k}^2\leq
  \sum_{k=0}^N\liminf_{n\to\infty}\|f_{nk}\|_{Z_k}^2\leq
  \liminf_{n\to\infty}\sum_{k=0}^N\|f_{nk}\|_{Z_k}^2\leq M^2.
$$
  Letting $N\to\infty$, we get the desired assertion. $\quad\Box$
\medskip

  We are now ready to give the proof of Theorem 1.1.
\medskip

  {\it Proof of Theorem 1.1}:\ \ We split the proof into four steps.
\medskip

  {\it Step 1}:\ \ We prove that if the assertion of Theorem 3.1 holds for
  $T=\delta$ for some small quantity $\delta>0$, then it also holds for
  any given $T>0$. This follows from a division of the time interval
  $[-T,T]$ and an induction argument. Indeed, let $m=T/\delta$ if $T/\delta$
  is an integer and $m=[T/\delta]+1$ otherwise. Let $I_j=[(j\!-\!1)\delta,
  j\delta]$, $j=1,2,\cdots, m\!-\!1$, $I_m=[(m\!-\!1)\delta,T]$, and $I_{-j}
  =-I_j$, $j=1,2,\cdots, m$. Since the length of each time interval $I_{\pm j}$
  is not larger than $\delta$, by assumption we see that the assertions $(i)$
  and $(ii)$ of Theorem 1.1 applies to each of these intervals provided
  $u_n(\cdot,t)$ weakly converges to $u(\cdot,t)$ in $L^2(\mathbb{R})$ for $t$
  equal to one of the two endpoints of this interval, but which follows from induction.
  Hence, the assertions $(i)$ and $(ii)$ of Theorem 1.1 holds for each of these
  intervals. Now, since for any $f\in L^{q'}([-T,T], L^{p'}(\mathbb{R}))$ ($p$,
  $q$ are as in $(i)$ of Theorem 1.1) we have
$$
  \int_{-T}^T\!\!\int_{-\infty}^\infty\!\! [u_n(x,t)-u(x,t)]f(x,t)dxdt
  =\sum_{|j|=1}^m\int_{I_j}\!\int_{-\infty}^\infty\!\!
  [u_n(x,t)-u(x,t)]f(x,t)dxdt,
$$
  the assertion $(i)$ follows immediately. Similarly, since for any
  $\varphi\in L^2(\mathbb{R})$,
$$
  \sup_{|t|\leq T}|(u_n(\cdot,t)-u(\cdot,t),\varphi)_{L^2}|=
  \sup_{1\leq |j|\leq m}\sup_{t\in I_j}|(u_n(\cdot,t)-u(\cdot,t),\varphi)_{L^2}|
$$
  the assertion $(ii)$ also follows immediately.
\medskip

  {\it Step 2}:\ \ By the result of Step 1 combined with a standard scaling
  argument, we see that we only need to prove Theorem 1.1 under the additional
  assumption that for $\epsln$ as in Theorem 2.3,
$$
  \|\phi_n\|_{L^2}\leq\epsln, \quad n=1,2,\cdots \quad
  \mbox{and} \quad  \|\phi\|_{L^2}\leq\epsln.
\eqno{(4.9)}
$$
  Thus, in what follows we always assume that this assumption is satisfied.
  Moreover, by density of $C([-T,T],C_0^\infty(\mathbb{R}))$ in $L^{q'}([-T,T],
  L^{p'}(\mathbb{R}))$ for any admissible pair $(p,q)$ and boundedness of the
  sequence $\{u_n\}$ in $L^q([-T,T],L^p(\mathbb{R}))$ (ensured by Corollary
  3.8), it can be easily seen that the assertion $(i)$ follows if we prove
  that
$$
  \mbox{for any}\;\; f\in C([-T,T],C_0^\infty(\mathbb{R})),\quad
  \displaystyle\lim_{n\to\infty}\int_{-T}^T\!\!\int_{-\infty}^\infty\!\!
  [u_n(x,t)-u(x,t)]f(x,t)dxdt=0.
\eqno{(4.10)}
$$
  Similarly, by density of $C_0^\infty(\mathbb{R})$ in $L^2(\mathbb{R}))$ and
  uniform boundedness of $\{u_n(\cdot,t)\}$ in $L^2(\mathbb{R})$ ensured by
  the $L^2$ conservation law, we see that the assertion $(ii)$ follows if we
  prove that
$$
  \mbox{for any}\;\; \varphi\in C_0^\infty(\mathbb{R}),\quad
  \displaystyle\lim_{n\to\infty}\sup_{|t|\leq T}|(u_n(\cdot,t)-u(\cdot,t),
  \varphi)_{L^2}|=0.
\eqno{(4.11)}
$$
\medskip

  {\it Step 3}:\ \ Due to (4.9), we have, by Theorem 2.3, the following
  expressions:
$$
  u_n=e^{-iU_{n0}}w_{n+}+e^{iU_{n0}}w_{n-}+w_{n0}+u_{n0},
  \quad n=1,2,\cdots,
\eqno{(4.12)}
$$
$$
  u=e^{-iU_0}w_+ + e^{iU_0}w_- + w_0 + u_0.
\eqno{(4.13)}
$$
  In what follows we prove that
$$
\left\{
\begin{array}{l}
  \mbox{for any}\;\; R>0\;\; \mbox{and}\;\; k,l\in\mathbb{Z}_+,\;\;
  \partial_t^k\partial_x^l u_{n0}\to\partial_t^k\partial_x^l u_0\\
  \mbox{uniformly on}\;\; [-R,R]\times [-T,T]\;\; \mbox{as}\;\;n\to\infty.
\end{array}
\right.
\eqno{(4.14)}
$$
  Note that if this assertion is proved, then it follows immediately that also
$$
\left\{
\begin{array}{l}
  \mbox{for any}\;\; R>0\;\; \mbox{and}\;\; k,l\in\mathbb{Z}_+,\;\;
  \partial_t^k\partial_x^l U_{n0}\to\partial_t^k\partial_x^l U_0\\
  \mbox{uniformly on}\;\; [-R,R]\times [-T,T]\;\; \mbox{as}\;\;n\to\infty.
\end{array}
\right.
\eqno{(4.15)}
$$

  Let $\phi_{{\rm low}}$ be as before, i.e., $\phi_{{\rm low}}=
  P_{{\rm low}}(\phi)$, and let $\phi_{n{\rm low}}=P_{{\rm low}}(\phi_n)$,
  $n=1,2,\cdots$. Then $u_{n0}$ and $u_0$ are respectively solutions of the
  following problems:
$$
\left\{
\begin{array}{l}
  \partial_t u_{n0}+\mathcal{H}\partial_x^2 u_{n0}+
  \partial_x(u_{n0}^2/2)=0, \quad
  x\in\mathbb{R},\;\;\; t\in\mathbb{R},\\
  u_{n0}(x,0)=\phi_{n{\rm low}}(x), \quad x\in\mathbb{R},
\end{array}
\right.
\eqno{(4.16)}
$$
$$
\left\{
\begin{array}{l}
  \partial_t u_0+\mathcal{H}\partial_x^2 u_0+
  \partial_x(u_0^2/2)=0, \quad
  x\in\mathbb{R},\;\;\; t\in\mathbb{R},\\
  u_0(x,0)=\phi_{{\rm low}}(x), \quad x\in\mathbb{R}.
\end{array}
\right.
\eqno{(4.17)}
$$
  Since the sequence $\{\phi_{n{\rm low}}\}$ is bounded in $L^2(\mathbb{R})$,
  by (2.2) we see that for any $k,l\in\mathbb{Z}_+$ there exists corresponding
  constant $C_{kl}(T)>0$ such that
$$
  \|\partial_t^k\partial_x^l u_{n0}\|_{L^2(\mathbb{R}\times [-T,T])}
  \leq\sqrt{2T}\sup_{|t|\leq T}\|\partial_t^k\partial_x^l
  u_{n0}(\cdot,t)\|_{L_x^2}\leq C_{kl}(T), \quad n=1,2,\cdots.
$$
  Hence, by using the compact embedding $H^{m+2}((-R,R)\times(-T,T))
  \hookrightarrow C^m([-R,R]\times [-T,T])$ for any $R>0$ and $m\in\mathbb{Z}_+$
  and a diagonalisation argument we see that there exists a subsequence
  $\{u_{n_k0}\}$ of $\{u_{n0}\}$, such that for any $R>0$ and $k,l\in
  \mathbb{Z}_+$, $\partial_t^k\partial_x^l u_{n_k0}$ is uniformly convergent in
  $[-R,R]\times [-T,T]$. Replacing $n$ with $n_k$ in (4.16) and then letting
  $k\to\infty$, we see that the limit function of $u_{n_k0}$
  is a smooth solution of the problem (4.17). Since $\phi_{{\rm low}}\in
  H^\infty(\mathbb{R})$, we know that the solution of (4.17) in $C([-T,T],
  H^\infty(\mathbb{R}))$ is unique, so that the limit function of $u_{n_k0}$ is
  $u_0$. Thus we have shown that
$$
\left\{
\begin{array}{l}
  \mbox{for any}\;\; R>0\;\; \mbox{and}\;\; k,l\in\mathbb{Z}_+,\;\;
  \partial_t^k\partial_x^l u_{n_k0}\to\partial_t^k\partial_x^l u_0\\
  \mbox{uniformly on}\;\; [-R,R]\times [-T,T]\;\; \mbox{as}\;\;n\to\infty.
\end{array}
\right.
$$
  Since the above argument is also valid when $\{u_{n0}\}$ is replaced by any
  of its subsequence, we see that the assertion (4.14) follows.
\medskip

  {\it Step 4}:\ \ By the assertions (4.14), (4.15) and the expressions (4.12),
  (4.13), it follows immediately that (4.10) and (4.11) will follow if we prove
  that
$$
  w_{n\alpha}\to w_{\alpha}\;\;\mbox{weakly in}\;\; L^q([-T,T],L^p(\mathbb{R}))
  \;\; \mbox{for any admissble pair}\;\; (p,q)
\eqno{(4.18)}
$$
  (in case either $p=\infty$ or $q=\infty$, weakly here refers to
  $\ast$-weakly), and
$$
  \mbox{for any}\;\;\varphi\in C_0^\infty(\mathbb{R}),\quad
  \displaystyle\lim_{n\to\infty}\sup_{|t|\leq T}
  |(w_{n\alpha}(\cdot,t)-w_{\alpha}(\cdot,t),\varphi)_{L^2}|=0,
\eqno{(4.19)}
$$
  where $\alpha=+,-,0$.

  To prove the assertion (4.18), we note that since $\phi_n\to\phi$ weakly in
  $L^2(\mathbb{R})$, using the assertion (4.15) we easily see that also
  $\psi_\alpha(\phi_n)\to\psi_\alpha(\phi)$ weakly in $L^2(\mathbb{R})$ for
  $\alpha=+,-,0$. Moreover, by Lemma 10.1 of \cite{IK} we see that
  $\{\psi_\alpha(\phi_n)\}$ ($\alpha=+,-,0$) are bounded in $\widetilde{H}^0$.
  The latter assertion implies that the sequences $\{w_{n\alpha}\}$ ($\alpha=+,
  -,0$) are bounded in $F_T^0$, which further implies, by Lemma 3.6, that for
  any admissible pair $(p,q)$ the sequences $\{w_{n\alpha}\}$ ($\alpha=+,-,0$)
  are also bounded in $L^q([-T,T],L^p(\mathbb{R}))$. If the assertion (4.18)
  does not hold for some admissible pair $(p,q)$ then it follows that there
  exist subsequences $\{w_{n_k\alpha}\}$ ($\alpha=+,-,0$) and functions
  $w_\alpha'\in L^q([-T,T],L^p(\mathbb{R}))$ ($\alpha=+,-,0$), $(w_{+}',w_{-}',
  w_{0}')\neq(w_{+},w_{-},w_{0})$, such that $w_{n_k\alpha}\to w_\alpha'$
  weakly in $L^q([-T,T],L^p(\mathbb{R}))$ (in case either $p=\infty$ or $q=
  \infty$ then weakly here refers to $\ast$-weakly, which will not be repeated
  later on). Since $\{w_{n_k\alpha}\}$ ($\alpha=+,-,0$) are bounded in
  $L^\infty([-T,T],L^2(\mathbb{R}))$ (by Lemma 3.6), by replacing them with
  subsequences of them when necessary, we may assume that also $w_{n_k\alpha}
  \to w_\alpha'$ ($\alpha=+,-,0$) weakly in $L^\infty([-T,T],L^2(\mathbb{R}))$,
  so that also $w_\alpha'\in L^\infty([-T,T],L^2(\mathbb{R}))$ ($\alpha=+,-,0$).
  Using Lemmas 3.2, 3.4 and a well-known compact embedding result, we easily
  deduce that there exists a subsequence, which we still denote as $\{w_{n_k
  \alpha}\}$, such that for any $R>0$, $w_{n_k\alpha}\to w_\alpha'$ strongly
  in $L^2([-R,R]\times [-T,T])$ ($\alpha=+,-,0$). Thus, by replacing
  $\{w_{n_k\alpha}\}$ with a subsequence when necessary, we may assume that
  $w_{n_k\alpha}\to w_\alpha'$ ($\alpha=+,-,0$) almost everywhere in $\mathbb{R}
  \times [-T,T]$. Now, From (2.4) we have
$$
  w_{n_k\alpha}(t)=W(t)\psi_\alpha(\phi_{n_k})+
  \int_0^t W(t-t')E_\alpha(w_{n_k+}(t'),w_{n_k-}(t'),w_{n_k0}(t'))dt',
  \quad  \alpha=+,-,0.
\eqno{(4.20)}
$$
  Here the second term on the right-hand sides of the above equations should
  be understood in the following sense: All partial derivatives in $x$ included
  in $E_\alpha$'s acting on any terms containing $w_{n_k +}$, $w_{n_k -}$ or
  $w_{n_k 0}$ should be either taken outside of the integral or moved to
  other terms containing only $u_{n_k0}$ and $U_{n_k0}$ by using integration by
  parts. For instance, recalling that (see (2.11) of \cite{IK})
$$
\begin{array}{rl}
  E_+(w_+,&w_-,w_0)=-e^{iU_0}P_{+{\rm high}}[\partial_x(e^{-iU_0}w_+
  +e^{iU_0}w_-+w_0)^2/2]
\\ [0.2cm]
  &-e^{iU_0}P_{+{\rm high}}\{\partial_x[u_0\cdot
  P_{-{\rm high}}(e^{iU_0}w_-)+u_0\cdot P_{{\rm low}}(w_0)]\}
\\ [0.2cm]
  &+e^{iU_0}(P_{-{\rm high}}+P_{{\rm low}})\{\partial_x[u_0\cdot
  P_{+{\rm high}}(e^{-iU_0}w_+)]\}
\\ [0.2cm]
  &+2iP_{-}\{\partial_x^2[e^{iU_0}P_{+{\rm high}}(e^{-iU_0}w_+)]\}
\\ [0.2cm]
  &-P_{+}(\partial_x u_0)\cdot w_+,
\end{array}
$$
  the equation in (4.20) for $\alpha=+$ should be understood to represent the
  following equation:
$$
\begin{array}{rcl}
  w_{n_k+}(t)&=&\displaystyle W(t)\psi_+(\phi_{n_k})-
  \partial_x\!\int_0^t W(t-t')e^{i U_{n_k 0}}P_{+{\rm high}}
  [(e^{-i U_{n_k 0}}w_{n_k +}+e^{i U_{n_k 0}}w_{n_k -}+w_{n_k 0})^2/2]dt'
\\ [0.5cm]
  &&\displaystyle+\!\int_0^t W(t-t') \partial_x(e^{i U_{n_k 0}})
  P_{+{\rm high}}[(e^{-i U_{n_k 0}}w_{n_k +}+e^{i U_{n_k 0}}w_{n_k -}
  +w_{n_k 0})^2/2]dt'
\\ [0.5cm]
  &&\displaystyle-\partial_x\!\int_0^t W(t-t')e^{i U_{n_k 0}}
  P_{+{\rm high}}[u_{n_k 0}\cdot P_{-{\rm high}}(e^{i U_{n_k 0}}w_{n_k -})+
  u_{n_k 0}\cdot P_{{\rm low}}(w_{n_k 0})]dt'
\\ [0.5cm]
  &&\displaystyle+\!\int_0^t W(t-t')\partial_x(e^{i U_{n_k 0}})
  P_{+{\rm high}}[u_{n_k 0}\cdot P_{-{\rm high}}(e^{i U_{n_k 0}}w_{n_k -})+
  u_{n_k 0}\cdot P_{{\rm low}}(w_{n_k 0})]dt'
\\ [0.5cm]
  &&\displaystyle+\partial_x\!\int_0^t W(t-t')e^{i U_{n_k 0}}
  (P_{-{\rm high}}+P_{{\rm low}})[u_{n_k 0}\cdot
  P_{+{\rm high}}(e^{-i U_{n_k 0}}w_{n_k +})]dt'
\\ [0.5cm]
  &&\displaystyle-\!\int_0^t W(t-t')\partial_x(e^{i U_{n_k 0}})
  (P_{-{\rm high}}+P_{{\rm low}})[u_{n_k 0}\cdot
  P_{+{\rm high}}(e^{-i U_{n_k 0}}w_{n_k +})]dt'
\\ [0.5cm]
  &&\displaystyle+2i\partial_x^2\!\int_0^t W(t-t')P_-[e^{i U_{n_k 0}}
  P_{+{\rm high}}(e^{-i U_{n_k 0}}w_{n_k +})]dt'
\\ [0.5cm]
  &&\displaystyle-\!\int_0^t W(t-t')P_+(\partial_x u_{n_k 0}\cdot w_{n_k +})dt';
\end{array}.
$$
  moreover, the equations in (4.20) for $\alpha=-$ and $\alpha=0$ should be
  understood similarly. Thus, letting $k\to\infty$ and using the Vitali convergence
  theorem (see Corollary A.2 in the appendix B), we see that $(w_{+}',w_{-}',
  w_{0}')$ satisfies the integral equations
$$
  w_{\alpha}'(t)=W(t)\psi_\alpha(\phi)+
  \int_0^t W(t-t')E_\alpha(w_{+}'(t'),w_{-}'(t'),w_{0}'(t'))dt',
  \quad  \alpha=+,-,0.
$$
  Note that these equations should be understood as (4.20) in the meaning
  explained above. Since by Lemma 4.1 we have $w_{\alpha}'\in F_T^0$
  ($\alpha=+,-,0$) and both $(w_{+}',w_{-}',w_{0}')$ and $(w_{+},w_{-},w_{0})$
  are in a small ball of $(F_T^0)^3$, by uniqueness of the solution of the above
  equation in a small ball of $(F_T^0)^3$ we conclude that $(w_{+}',w_{-}',
  w_{0}')=(w_{+},w_{-},w_{0})$, which is a contradiction.

  The argument for the proof of (4.19) is similar. Indeed, let $v_{n_k\alpha}^1$
  and $v_{n_k\alpha}^2$ denote the first and the second terms on the right-hand
  side of (4.20), respectively, and by $v^1$ and $v^2$ the corresponding terms
  in $w_\alpha'$. It can be easily seen that
$$
  \lim_{k\to\infty}\sup_{|t|\leq T}|(v_{n_k\alpha}^1(\cdot,t),\varphi)_{L_x^2}
  -v^1(\cdot,t),\varphi)_{L_x^2}|=0.
\eqno{(4.21)}
$$
  To treat $(v_{n_k\alpha}^2(\cdot,t),\varphi)_{L_x^2}$ we only need to move all
  partial derivatives in $x$ contained in $E_\alpha$'s either to terms expressed
  in $u_{n_k0}$ and $U_{n_k0}$ by using integration by parts, or to the test
  function $\varphi$, also by using integration by parts. With this trick in
  mind, we can also prove that
$$
  \lim_{k\to\infty}\sup_{|t|\leq T}|(v_{n_k\alpha}^2(\cdot,t),\varphi)_{L_x^2}
  -(v^2(\cdot,t),\varphi)_{L_x^2}|=0.
\eqno{(4.22)}
$$
  We omit the details. combining (4.21) and (4.22), we see that the assertion
  (4.19) follows. This completes the proof of Theorem 1.1. $\quad\Box$
\medskip

  {\bf Remark}\ \  For the modified Benjamin-Ono equation (1.3), it has been
  proved by Kenig and Takaoka in \cite{KT} that its initial value problem is
  globally well-posed in the Sobolev space $H^{1/2}(\mathbb{R})$, whereas the
  solution operator of a such problem is not uniformly continuous in any Sobolev
  spaces $H^{s}(\mathbb{R})$ of index $s<1/2$ (so that $H^{1/2}(\mathbb{R})$ is
  a borderline space for the local well-posedness theory of this equation). It
  is thus natural to ask if the flow map of this equation in
  $H^{1/2}(\mathbb{R})$ is weakly continuous. The answer to this question is
  affirmative. The proof is as follows: Let $\phi_n$ ($n=1,2,\cdots$) be a
  sequence of functions in $H^{1/2}(\mathbb{R})$ which is weakly convergent,
  and let $\phi$ be its limit. Let $u_n$ and $u$ be the solutions of the
  equation (1.3) in $C(\mathbb{R},H^{1/2}(\mathbb{R}))$ such that $u_n|_{t=0}=
  \phi_n$ ($n=1,2,\cdots$) and $u|_{t=0}=\phi$. Then for any $T>0$, $\{u_n\}$
  is bounded in $L^\infty([-T,T],H^{1/2}(\mathbb{R}))$. Using the equation (1.3),
  we then deduce that $\{\partial_t u_n\}$ is bounded in $L^\infty([-T,T],
  H^{-3/2}(\mathbb{R}))$. It follows that there exists a subsequence
  $\{u_{n_k}\}$ such that for any $R>0$, $\{u_{n_k}\}$ is strongly convergent
  in $L^2([-R,R]\times[-T,T])$ and, consequently, by replacing $\{u_{n_k}\}$
  with a suitable subsequence of it, we may assume that $\{u_{n_k}\}$
  converges almost everywhere in $\mathbb{R}\times [-T,T]$. Thus by following
  the approach developed in \cite{GoM} we obtain the desired assertion. (One
  needs to, in addition, observe that the uniqueness in \cite{KT} easily
  extends to solutions of the integral equation in $C([-T,T],H^{1/2}(\mathbb{R}))
  \cap X^{1/2}$, where $X^{1/2}$ is the space in \cite{KT}). We are grateful to
  one of the anonymous referees for pointing to us this proof.
\medskip

\section*{Appendix: Vitali convergence theorem}
\renewcommand{\theequation}{$B$.\arabic{equation}}

\hspace*{2em}
  {\bf Theorem A.1} (Vitali convergence theorem, cf. \cite{F})\ \ {\em Let $X$
  be a measurable set. Let $u_n\in L^1(X)$, $n=1,2,\cdots$. Assume that the
  following three conditions are satisfied:\\
\hspace*{2em}  $(a)$ $u_n$ converges to $u$ in measure.\\
\hspace*{2em}  $(b)$ For any $\epsln>0$ there exists corresponding $M>0$ such that
$$
  \int_{\{|u_n(x)|>M\}}|u_n(x)|dx<\epsln \quad \mbox{for all}\;\; n\in\mathbb{N}.
$$
\hspace*{2em}  $(c)$ For any $\epsln>0$ there exists corresponding measurable
  subset $E$ of $X$ with $\meas(E)<\infty$,\\
\hspace*{2em} such that
$$
  \int_{X\backslash E}|u_n(x)|dx<\epsln \quad \mbox{for all}\;\; n\in\mathbb{N}.
$$
  Then $u\in L^1(X)$ and
$$
  \lim_{n\to\infty}\int_E u_n(x)dx=\int_E u(x)dx.
$$
}

  {\bf Remark} \ \ If $\meas(X)<\infty$, then the condition $(c)$ is clearly
  satisfied by any sequence of measurable functions on $X$: We may choose $E=X$.

  What we used in the proof of Theorem 1.1 is the following corollary of the
  above theorem:
\medskip

  {\bf Corollary A.2}\ \ {\em Let $E$ be a measurable set, $\meas(E)<\infty$. Let
  $1<p<\infty$ and $u_n\in L^p(E)$, $n=1,2,\cdots$. Assume that $(i)$ $u_n$
  converges to $u$ in measure, and $(ii)$ $\{u_n\}$ is bounded in $L^p(E)$.
  Then $u\in L^p(E)$, and for any $1\leq q<p$ we have
$$
  \lim_{n\to\infty}\|u_n-u\|_q=0.
\eqno{(A.1)}
$$
}

  {\em Proof}:\ \ The assertion that $u\in L^p(E)$ follows from Fatou's lemma.
  To prove $(A.1)$ we assume that $\|u_n\|_p\leq C$ for all $n\in\mathbb{N}$.
  Then we also have $\|u\|_p\leq C$, by Fatou's lemma. Thus, for any $M>0$ we
  have
$$
  \int_{\{|u_n(x)\!-\!u(x)|>M\}}\!|u_n(x)\!-\!u(x)|^q dx\leq M^{-(p-q)}\!
  \int_{E}\!|u_n(x)\!-\!u(x)|^p dx\leq (2C)^p M^{-(p-q)},
$$
  which implies that
$$
  \lim_{M\to\infty}\sup_{n\in\mathbb{N}}\int_{\{|u_n(x)\!-\!u(x)|>M\}}\!
  |u_n(x)\!-\!u(x)|^q dx=0.
$$
  Hence, the desired assertion follows from Theorem A.1. $\quad\Box$


{\small
\begin{thebibliography}{99}
\bibitem[1]{B} T. B. Benjamin, Internal waves of permanent form in fluids of
  great depth, \textit{J. Fluid Mech.}, \textbf{29}(1967), 559--592.
\bibitem[2]{BeL} J. Bergh and J. L\"{o}fstr\"{o}m, \textit{Interpolation
  Spaces: An Introduction}, Springer-Verlag, New York, 1976.
\bibitem[3]{BiL} H. A. Biagioni and F. Linares, Ill-posedness for the
  derivative Schr\"{o}dinger and generalized Benjamin-Ono equations,
  \textit{Trans. Amer. Math. Soc.}, \textbf{353}(2001), 3649--3659.
\bibitem[4]{B} J. Bourgain, \textit{Global Solutions of Nonlinear
  Schr\"{o}dinger Equations}, Amer. Math. Soc. Colloquium Publ.,
  vol. 46, Providence, Rhode Island, 1999.
\bibitem[5]{BP} N. Burq and F. Planchon, On well-posedness for the
  Benjamin-Ono equation, \textit{Math. Annalen}, \textbf{340}(2008),
  497--542
\bibitem[6]{CCT} M. Christ, J. Colliander and T. Tao, Asymptotics,
  frequency modulation and low-regularity ill-posedness for canonical
  defocusing equations, \textit{Amer. J. Math.}, \textbf{125}(2003),
  1235--1293.
\bibitem[7]{CKS} J. Colliander, C. E. Kenig and G. Staffilani, Local
  well-posedness for dispersion-generalized Benjamin-Ono equations,
 \textit{Differential Integral Equations}, \textbf{16}(2003), 1441--1472.
\bibitem[8]{GV1} J. Ginibre and G. Velo, Smoothing properties and
  existence of solutions for the generalized Benjamin-Ono equation,
 \textit{J. Diff. Equations}, \textbf{93}(1991), 150--212.
\bibitem[9]{GV2} J. Ginibre and G. Velo, Properties de lissage et
  existence de solutions pour l¡¯equation de Benjamin-Ono generalisee,
 \textit{C. R. Acad. Sci. Paris Ser. I. Math.}, \textbf{308}(1989), 309--314.
\bibitem[10]{GV3} J. Ginibre and G. Velo, Commutator expansions and
  smoothing properties of generalized Benjamin-Ono equations,
 \textit{Ann. Inst. H. Poincare, Phys. Theor.}, \textbf{51}(1989), 221--229.
\bibitem[11]{GoM} O. Goubet and L. Molinet, Global weak attractor for weakly
  damped nonlinear Schr\"{o}dinger equations in $L^2(R)$, \textit{Nonlinear
  Anal.}, \textbf{71}(2009), 317--320.
\bibitem[12]{F} G. B. Folland, \textit{Real Analysis}, 2nd ed., New York:
  John Wiley \& Sons Inc., 1999.
\bibitem[13]{IK} A. D. Ionescu and C. E. Kenig, Global well-posedness
  of the Benjamin-Ono equation in low-regularity spaces, \textit{J. Amer.
  Math. Soc.}, \textbf{20}(2007), 753--798.
\bibitem[14]{I} R. J. Iorio, On the Cauchy problem for the Benjamin-Ono
  equation, \textit{Comm. Partial Differential Equations}, \textbf{11}(1986),
  1031--1081.
\bibitem[15]{KPV2} C. E. Kenig, G. Ponce and L. Vego, Well-posedness and
  scattering results for the generalized Korteweg-de Vries equation via
  contraction principle, \textit{Comm. Pure Appl. Math.}, \textbf{46}(1993),
  527--620.
\bibitem[16]{KPV3} C. E. Kenig, G. Ponce and L. Vega, On the
  generalized Benjamin-Ono equation, \textit{Trans. Amer. Math. Soc.},
  \textbf{342}(1994), 155--172.
\bibitem[17]{KPV4} C. E. Kenig, G. Ponce and L. Vega, On the ill-posedness
  of some canonical  dispersive equations, \textit{Duke Math. J.},
  \textbf{106}(2001), 617--633.
\bibitem[18]{KT} C. E. Kenig and H. Takaoka, Global well-posedness of the
  modified Benjamin-Ono equation with initial data in $H^{1/2}$,
  \textit{Int. Math. Res. Notes}, \text{2006}(2006), 1--44.
\bibitem[19]{KK} C. E. Kenig and K. D. Koenig, On the local well-posedness
  of the Benjamin-Ono and modified Benjamin-Ono equations, \textit{Math. Res.
  Lett.}, \textbf{10}(2003), 879--895.
\bibitem[20]{KM} C. E. Kenig and Y. Martel, Asymptotic stability of
  solitons for the Benjamin-Ono equation, \textit{Revista Mat.
  Iberoamericana}, to appear (see also arXiv:0803.3683).
\bibitem[21]{KMR} C. E. Kenig, Y. Martel and L. Robbiano, Well-posedness and
  blow-up in the energy space for critical dispersive generalized Benjamin-Ono
  equations, in preparation.
\bibitem[22]{KoTz}  H. Koch and N. Tzvetkov, On the local well-posedness of
  the Benjamin-Ono equation in $H^s(R)$, \textit{Int. Math. Res. Not.},
  2003, no. 26, 1449--1464.
\bibitem[23]{MM1} Y. Martel and F. Merle, A Liouville theorem for the
  critical generalized Korteweg-de Vries equation, \textit{J. Math. Pures
  Appl.}, \textbf{79}(2000), 339--425.
\bibitem[24]{MM2} Y. Martel and F. Merle, Instability of solitons for the
  critical generalized Korteweg-de Vries equation, \textit{Geom. Funct.
  Anal.}, \textbf{38}(2001), 759--781.
\bibitem[25]{MM3} Y. Martel and F. Merle, Asymptotic stability of solitons
  for subcritical generalized KdV equations, \textit{Arch. Rat. Mech. Anal.},
  \textbf{157}(2001), 219--254.
\bibitem[26]{MM4} Y. Martel and F. Merle, Blow up in finite time and
  dynamics of blow up solutions for the $L^2$-critical generalized KdV equation,
  \textit{J. Amer. Math. Soc.}, \textit{15}(2002), 617--664.
\bibitem[27]{M1} L. Molinet, Global well-posedness in the energy space
  for the Benjamin-Ono equation on the circle, \textit{Math. Ann.},
 \textbf{337}(2007), 353--383.
\bibitem[28]{M2} L. Molinet, Global well-posedness in $L2$ for the periodic
  Benjamin-Ono equation, \textit{Amer. J. Math.}, \textbf{130}(2008), 635--683.
\bibitem[29]{M3} L. Molinet, On ill-posedness for the one-dimensional
  periodic cubic Schrodinger equation, to appear in \textit{Math. Res. Let.}
  \textbf{16}(2009), 111--120
\bibitem[30]{M4} L. Molinet, Sharp ill-posedness result for the periodic
  Benjamin-Ono equation, arXiv: 0811.0505.
\bibitem[31]{MR1} L. Molinet and F. Ribaud, Well-posedness results for the
  generalized Benjamin-Ono equation with small initial data, \textit{J. Math.
  Pures Appl.}, \textbf{83}(2004), 277--311.
\bibitem[32]{MR2} L. Molinet and F. Ribaud, Well-posedness results for the
  generalized Benjamin-Ono equation with arbitrary large initial data,
 \textit{Int. Math. Res. Not.}, 2004, no. 70, 3757--3795.
\bibitem[33]{Ono} H. Ono, Algebraic solitary waves in stratified fluids,
  \textit{J. Phys. Soc. Japan}, \textbf{39}(1975), 1082--1091.
\bibitem[34]{P} G. Ponce, On the global well-posedness of the Benjamin-Ono equation,
  \textit{Differential Integral Equations}, \textbf{4}(1991), 527--542.
\bibitem[35]{T} T. Tao, Global well-posedness of the Benjamin-Ono
  equation in $H^1(\mathbb{R})$, \textit{J. Hyperbolic Diff. Equa.},
  \textbf{1}(2004), 27--49.
\bibitem[36]{Ts} Y. Tsutsumi, $L^2$-solutions for nonlinear Schr\"{o}dinger
  equations and nonlinear groups, \textit{Funk. Ekva.}, \textbf{30}(1987),
  115--125.
\end{thebibliography}
}
\end{document}